\documentclass{conm-p-l}

\usepackage{amsfonts}
\usepackage{amsmath}
\usepackage{amssymb}
\usepackage{amscd}
\usepackage{mathrsfs}  % This package allows the use of script letter. Type \mathscr to invoke math script.
\usepackage{amsbsy}
\usepackage{amsthm}
\usepackage{graphicx}
\usepackage{fancyhdr}
\usepackage{color}

\newtheorem{theorem}{Theorem}[section]
\newtheorem{lemma}[theorem]{Lemma}
\newtheorem{corollary}[theorem]{Corollary}

\theoremstyle{definition}
\newtheorem{definition}[theorem]{Definition}
\newtheorem{example}[theorem]{Example}
\newtheorem{proposition}[theorem]{Proposition}

\theoremstyle{remark}
\newtheorem{remark}[theorem]{Remark}

\numberwithin{equation}{section}

\newcommand{\be}{\begin{equation}}
\newcommand{\ee}{\end{equation}}
\newcommand{\inn}[2]{{\langle #1,#2 \rangle}}
\newcommand{\st}{\,|\,}
\newcommand{\R}{\mathbb{R}}
\newcommand{\N}{\mathbb{N}}
\newcommand{\Z}{\mathbb{Z}}
\newcommand{\wW}{\widetilde{\mathcal{W}}}

\begin{document}

\title{Coxeter Groups and Wavelet Sets}

\author{David R. Larson}
\address{Department of Mathematics, Texas A \& M University, College Station, Texas 77843}
%\curraddr{Department of Mathematics and Statistics,
%Case Western Reserve University, Cleveland, Ohio 43403}
\email{larson@math.tamu.edu}

\author{Peter Massopust}
\address{GSF - National Research Center for Environment and Health, Institute of Biomathematics and Biometry, and Centre of Mathematics M6,
Technische Universit\"{a}t M\"{u}nchen, Germany}
\email{massopust@ma.tum.de}
\thanks{The first author was partially supported by grants from the National Science Foundation and the second author was partially supported by the
grant MEXT-CT-2004-013477, Acronym MAMEBIA, of the European Commission. Both authors participated in NSF supported Workshops in Linear Analysis and
Probability, Texas A\&M University.}

%    General info
\date{January 1, 1994 and, in revised form, June 22, 1994.}



\begin{abstract}
A traditional wavelet is a special case of a vector in a separable Hilbert space that generates a basis under the action of a  system of unitary
operators defined in terms of translation and dilation operations. A Coxeter/fractal-surface wavelet is obtained by defining fractal surfaces on
foldable figures, which tesselate the embedding space by reflections in their bounding hyperplanes instead of by translations along a lattice.
Although both theories look different at their onset, there exist connections and communalities which are exhibited in this semi-expository paper.
In particular, there is a natural notion of a dilation-reflection wavelet set.  We prove that dilation-reflection wavelet sets exist for arbitrary
expansive matrix dilations, paralleling the traditional dilation-translation wavelet theory. There are certain measurable sets which can serve
simultaneously as dilation-translation wavelet sets and dilation-reflection wavelet sets, although the orthonormal structures generated in the two
theories are considerably different.
\end{abstract}
\maketitle
\section{Introduction}
This article is meant to mesh together two distinct approaches to wavelet theory: the traditional dilation-translation approach that was the subject
of the memoir \cite{DL} co-authored by the first author, and the Coxeter/fractal-surface approach that was developed by the second author in
\cite{GHM1,GHM2} and the book \cite{M2}.  The two approaches seem distinctly different under initial scrutiny.  However, it turns out that both
approaches carry a natural notion of ``wavelet set", and while these have different meanings in the two theories, we have recently discovered, with
some surprise, that there is more than a little connection between the two notions of ``wavelet set."  Indeed, in the plane $\mathbb{R}^2$, with
dilation (scale factor) $2$, the ``dyadic case", some of the known wavelet sets in the dilation-translation theory, and in particular the ``wedding
cake set" (c.f. \cite{DL}, p. 59, and \cite{DLS2}) and the ``four-corners set" (c.f. \cite{DL}, p.57, and  \cite{DLS2}) are both also wavelet sets in
the Coxeter/fractal-surface multiresolution-analysis theory. We remark that Coxeter/fractal-surface ``wavelet sets", as such, were not formally
defined  in the book \cite{M2} and were not part of the original Coxeter-MRA theory that was formally defined and developed in \cite{GHM1, GHM2} and
\cite{M2}.  However, in a project we began three years ago, while both of us were participants at an international conference on Abstract and
Applied Analysis in Hanoi, Vietnam, we decided to pursue the possibility of a connection between the subjects of our two talks. This resulted in
developing the \cite{GHM1, GHM2, M2} theory a bit further, including a proper notion of ``wavelet set" in that context.
\par
The structure of this paper is as follows. In Section 2, we review some aspects of traditional wavelet theory, introduce the concept of wavelet set,
translation and dilation congruence, and abstract dilation-translation pair, and summarize some of the results from the theory of
dilation-translation wavelet sets. A class of fractal functions and fractal surfaces is introduced in Section 3 and some of their properties are
mentioned. Section 4 deals with Coxeter groups, affine Weyl groups, and foldable figures to the extend that is necessary for this paper. Fractal
surfaces on foldable figures are defined in Section 5 and are shown to generate multiresolution analyses in Section 6. The concept of
dilation-reflection wavelet set is introduced in Section 7 and it is shown that these wavelet sets exist for arbitrary expansive matrix dilations.
We exhibit two examples of measurable sets that serve simultaneously as dilation-translation as well as dilation-reflection wavelet sets, although
the orthonormal structures they define are intrinsically different. In Section 8, we consider some questions of a general nature regarding
dilation-reflection wavelet sets and pose two open problems.
\section{Some Aspects of Traditional Wavelet Theory}
A traditional wavelet is a special case of a vector in a separable Hilbert space that generates a basis under the action of a collection, or
``system", of unitary operators defined in terms of translation and dilation operations. A traditional wavelet set is a measurable set whose
characteristic function, scaled appropriately, is the Fourier Transform of a (single) orthonormal wavelet.  Multi-versions have been studied, as
well as frame analogues.  In \cite{DLS1}, Dai and Speegle, together with the first author, found a proof of the existance of wavelet sets in the
plane and in higher dimensions.  The announcement of the existance  of such wavelet sets stimulated the construction of examples of such sets by
several authors, beginning with Soardi and Weiland \cite{SW}, and then Baggett, Medina and Merrill \cite{BMM}, and Benedetto and Leon \cite{BL}.
Wavelet sets have been a part of the wavelet literature for a number of years now.
\par
A \textit{dyadic orthonormal} wavelet in one dimension is a unit vector $\psi \in L^2(\mathbb{R}, \mu)$, with $\mu$ Lebesgue measure, with the property that the set
\begin{equation}
\{2^{\frac{n}{2}}\psi(2^n t - \ell) \st n,\ell \in \mathbb{Z}\}
\end{equation}
of all integral translates of $\psi$ followed by dilations by arbitrary integral powers of $2$, is an orthonormal basis for $L^2(\mathbb{R},\mu)$. The term \emph{dyadic} refers to the dilation factor ``$2$".  The term \emph{mother wavelet} is also used in the literature for $\psi$. Then the functions $$\psi_{n,\ell} (t):= 2^{\frac{n}{2}} \psi(2^n t - \ell)$$ are called elements of the wavelet basis generated by the ``mother".  The functions $\psi_{n,\ell}$ will not themselves be mother wavelets unless $n = 0$. Let $T$ and $D$ be the translation (by $1$) and dilation (by $2$) unitary operators in $B(L^2(\mathbb{R}))$, the Banach space of bounded linear operators from $L^2(\mathbb{R})$ to itself, given by $(Tf)(t) = f(t-1)$ and $(Df)(t) = \sqrt{2}f(2t)$.  Then
$$
2^{\frac{n}{2}}\psi(2^nt - \ell) = (D^nT^\ell \psi)(t)
$$
for all $n,\ell \in \mathbb{Z}$.  Operator-theoretically, the operators $T, D$ are \textit{bilateral shifts} of \textit{infinite multiplicity}.  It is obvious that $L^2([0,1])$, considered as a subspace of $L^2(\mathbb{R})$, is a complete wandering subspace for  $T$, and that  $L^2([-2, -1] \cup [1, 2])$ is a complete wandering subspace for $D$.
\par
A {\em complete wandering subspace} $E$ for a unitary operator $D$ acting on a Hilbert space $H$ is a closed subspace of $H$ for which $\{D^n E \st  n \in \mathbb{Z}\}$ is an orthogonal decomposition of $H$.  A unitary operator is called a {\em bilateral shift} if it has a complete wandering subspace.  In that case its multiplicity is defined to be the dimension of the wandering subspace.
\par
An abstract interpretation is that, since $D$ is a bilateral shift
it has (many) complete wandering subspaces,  and a wavelet for the
system is a vector $\psi$ whose translation space (that is, the
closed linear span of $\{T^k\st k \in \mathbb{Z}\}$ is a complete
wandering subspace for $D$.  Hence $\psi$ must generate an
orthonormal basis for the entire Hilbert space under the action of
the unitary system.
\par
In one dimension, there are non-dyadic orthonormal wavelets: i.e.
wavelets for all possible dilation factors besides $2$ (the dyadic
case). We said ``possible", because the scales $\{0,1,-1\}$ are
excluded as scales because the dilation operators they would
introduce are not bilateral shifts. All other real numbers for
scales yield wavelet theories. In \cite{DL}, Example 4.5 (x), a family of
examples is given of three-interval wavelet sets (and hence
wavelets) for all scales $d \geq 2$, and it was noted there that
such a family also exists for dilation factors $1 < d \leq 2$.
\par
Let $1  \leq n < \infty$, and let $A$ be an $n\times n$  real invertible matrix.   The most tractable such matrices for dilations in wavelet theory are those that are {\em expansive}. 
\begin{remark}
There are at least six equivalent characterizations of the property ``expansive" for an $n \times n$ real invertible matrix $A$.  It may be good to give them here for sake of exposition (and it can make a good student exercise to verify these equivalences).  Any of the first five can be (and sometimes is) taken as the definition, and the sixth is particularly useful for wavelet theory.  The first characterization is that all (complex) eigenvalues of $A$ have modulus $ > 1$.  A second is that $\bigcup \{A^\ell B_1 \st \ell \in \mathbb{N}\} = \mathbb{R}^n$, where $B_1 := B_1(0)$ is the open unit ball of  $\mathbb{R}^n$. A third  is that for each nonzero vector $x \in \mathbb{R}^n$, the sequence of norms $\{\|A^\ell x\| \st \ell \in \mathbb{N}\}$ is unbounded. A fourth characterization is that the sequence of norms$\{||A^{-\ell}|| : \ell \in\mathbb{N}\}$ converges to 0. A fifth (which is a quantitative version of the second) is that for each neighborhood $N$ of $0$ and each $r > 0$, there exists an $\ell \in \mathbb{N}$ such that $B_0(r) \subseteq  A^\ell N$.  And the sixth is that for each open set $F$ which is bounded away from $0$, and each $r > 0$, there exists an $\ell \in \mathbb{N}$ such that  $A^\ell F$  contains a ball $B_r(p)$ of radius $r$ and some center $p$ (where $p$ depends on $F$ and $r$).
\end{remark}
By a {\em dilation - A regular--translation orthonormal wavelet} we mean a function $\psi \in L^2(\mathbb{R}^n)$ such that
\begin{equation}
\{|\det(A)|^{\frac{n}{2}} \psi(A^n t - \ell) \st n\in \mathbb{Z}, \ell\in\mathbb{Z}^n, i = 1,\ldots, n\}
\end{equation}
where $\ell = (\ell_1, \ell_2, ..., \ell_n)^\top$, is an orthonormal basis for $L^2(\mathbb{R}^n ; m)$.  (Here $m$ is product Lebesgue measure, and the superscript ${}^\top$ means transpose.)If $A \in M_n(\mathbb{R})$ is invertible (so in particular if $A$ is expansive), then the operator defined by
\begin{equation}\label{e1}
(D_Af)(t) = |\det A|^{\frac12} f(At)
\end{equation}
for $f \in L^2(\mathbb{R}^n)$, $t \in \mathbb{R}^n$, is \emph{unitary}. For $1 \leq i \leq n$, let $T_i$ be the unitary
operator determined by translation by $1$ in the $i^{th}$ coordinate direction.  The set (5) above is then
\begin{equation}
\{D^k_A T^{\ell_1}_1 \cdot\cdot\cdot T^{\ell_n}_n \psi \st k,\ell_i \in \mathbb{Z}\}
\end{equation}
\par
If the dilation matrix $A$ is expansive, but the translations are
along some oblique lattice, then there is an invertible real $n
\times n$ matrix $T$ such that conjugation with $D_T$ takes the
entire wavelet system to a regular-translation expansive-dilation
matrix. This is easily worked out, and was shown in detail in \cite{ILP98} in the context of working out a complete theory of unitary
equivalence of wavelet systems. Hence the wavelet theories are
equivalent.
\par
Much work has also been accomplished concerning the existence of
wavelets for dilation matrices $A$ which are not expansive.  But
there is no need to go into that here.
\subsection{Fourier Transform}
We will use the following form of the Fourier--Plancherel transform
$\mathscr{F}$ on $\mathcal{H} = L^2(\mathbb{R})$, which is a form
that is \emph{normalized} so it is a unitary transformation, a
property that is desirable for our treatment. If $f,g\in
L^1(\mathbb{R}) \cap L^2(\mathbb{R})$ then
\begin{equation}\label{eq18}
(\mathscr{F}f)(s) := \frac1{\sqrt{2\pi}} \int_{\mathbb{R}} e^{-ist}
f(t)dt := \hat f(s),
\end{equation}
and
\begin{equation}\label{eq19}
(\mathscr{F}^{-1}g)(t) = \frac1{\sqrt{2\pi}} \int_{\mathbb{R}}
e^{ist}g(s)ds.
\end{equation}
We have
\[
(\mathscr{F} T_\alpha f)(s) = \frac1{\sqrt{2\pi}} \int_{\mathbb{R}}
e^{-ist} f(t-\alpha)dt = e^{-is\alpha} (\mathscr{F}f)(s).
\]
So $\mathscr{F} T_\alpha \mathscr{F}^{-1} g=e^{-is\alpha}g$. For
$A\in {B}(\mathcal{H})$ let $\hat A$ denote
$\mathscr{F}A\mathscr{F}^{-1}$. Thus
\begin{equation}\label{eq20}
\widehat T_\alpha = M_{e^{-i\alpha s}},
\end{equation}
where for $h\in L^\infty$ we use $M_h$ to denote the multiplication
operator $f\mapsto hf$. Since $\{M_{e^{-i\alpha s}}\st \
\alpha\in\mathbb{R}\}$ generates the m.a.s.a. (maximal abelian self adjoint operator algebra)\
$\mathcal{D}(\mathbb{R}) := \{M_h\st \ h\in
L^\infty(\mathbb{R})\}$ as a von Neumann algebra.  Let $\mathcal{A}_T$ denote the von Neumann algebra generated by $\{T_\alpha : \alpha\in
\mathbb{R}\}$. We then have
\[
\mathscr{F}\mathcal{A}_T \mathscr{F}^{-1} = \mathcal{D}(\mathbb{R}).
\]
Similarly,
\begin{align*}
(\mathscr{F}D^nf)(s) &= \frac1{\sqrt{2\pi}} \int_{\mathbb{R}}
e^{-ist} (\sqrt
2)^n f(2^nt)dt\\
&= (\sqrt 2)^{-n}\cdot \frac1{\sqrt{2\pi}} \int_{\mathbb{R}}
e^{-i2^{-n}st}
f(t)dt\\
&= (\sqrt 2)^{-2} (\mathscr{F}f)(2^{2^{-n}s}) =
(D^{-n}\mathscr{F}f)(s).
\end{align*}
So $\widehat D^n = D^{-n} = D^{*n}$. Therefore,
\begin{equation}\label{eq21}
\widehat D = D^{-1} = D^*.
\end{equation}
\par
Wavelet sets belong to the theory of wavelets via the Fourier
transform. As mentioned earlier, we define a {\em wavelet set\/ in $\mathbb{R}$} to
be a measurable subset $E$ of $\mathbb{R}$ for which
$\frac1{\sqrt{2\pi}} \chi_E$ is the Fourier transform of a wavelet.
The wavelet $\widehat\psi_E := \frac1{\sqrt{2\pi}}\chi_E$ is called
$s$-{\em elementary\/} in \cite{DL}. The class of wavelet sets was also
discovered and systematically explored completely independently, and
in about the same time period, by Guido Weiss (Washington
University), his colleague and former student E. Hernandez (U. Madrid), and his students X. Fang and X. Wang
\cite{HWW,FW}.  In this theory the
corresonding wavelets are are called MSF (minimally supported
frequency) wavelets.
\subsection{Shannon Wavelet}
The two most elementary dyadic orthonormal wave\-lets are the well-known {\em Haar wavelet\/} and {\em Shannon's wavelet\/} (also called the Little\-wood--Paley wavelet).  The Haar wavelet is the prototype of a large class of wavelets, and is the function given in the time domain by  $\psi =  \chi_{[0,1/2)} - \chi_{[1/2, 1)}$. The Shannon set (\ref{eq33}) is the prototype of the class of wavelet sets.
\par
Shannon's wavelet is the $L^2(\mathbb{R})$-function $\psi_S$ with Fourier
transform $\widehat\psi_S = \frac1{\sqrt{2\pi}} \chi_{E_0}$ where
\begin{equation}\label{eq33}
E_0 = [-2\pi, -\pi) \cup [\pi,2\pi).
\end{equation}
The argument that $\widehat\psi_S$ is a wavelet is in a way even
more transparent than for the Haar wavelet. And it has the advantage
of generalizing nicely. For a simple argument, start from the fact
that the set of exponentials
\[
\{e^{i\ell s}\st \ell\in \mathbb{Z}\}
\]
restricted to $[0,2\pi]$ and normalized by $\frac1{\sqrt{2\pi}}$ is
an orthonormal basis for $L^2[0,2\pi]$. Write $E_0 = E_-\cup E_+$
where $E_- = [-2\pi, -\pi)$, $E_+ = [\pi,2\pi)$. Since $\{E_- +2\pi,
E_+\}$ is a partition of $[0,2\pi)$ and since the exponentials
$e^{i\ell s}$ are invariant under translation by $2\pi$, it follows
that
\begin{equation}\label{eq34}
\left\{\frac{e^{i\ell s}}{\sqrt{2\pi}}\Big|_{E_0}\st \ell\in
\mathbb{Z}\right\}
\end{equation}
is an orthonormal basis for $L^2(E_0)$. Since $\widehat T =
M_{e^{-is}}$, this set can be written
\begin{equation}\label{eq35}
\{\widehat T^\ell \widehat\psi_s\st \ \ell\in \mathbb{Z}\}.
\end{equation}
Next, note that any ``dyadic interval'' of the form $J = [b,2b)$,
for some $b>0$ has the property that $\{2^nJ\st \
n\in\mathbb{Z}\}$, is a partition of $(0,\infty)$. Similarly, any
set of the form
\begin{equation}\label{eq36}
\mathcal{K} = [-2a,-a)\cup [b,2b)
\end{equation}
for $a,b>0$, has the property that
\[
\{2^n\mathcal{K}\st \ n\in \mathbb{Z}\}
\]
is a partition of $\mathbb{R}\backslash\{0\}$. It follows that the
space $L^2(\mathcal{K})$, considered as a subspace of
$L^2(\mathbb{R})$, is a complete wandering subspace for the dilation
unitary $(Df)(s) = \sqrt 2\ f(2s)$. For each $n\in \mathbb{Z}$,
\begin{equation}\label{eq37}
D^n(L^2(\mathcal{K})) = L^2(2^{-n}\mathcal{K}).
\end{equation}
So $\displaystyle{\bigoplus_{n\in\Z}} D^n(L^2(\mathcal{K}))$ is a direct sum decomposition
of $L^2(\mathbb{R})$. In particular $E_0$ has this property. So
\begin{equation}\label{eq38}
D^n\left\{\frac{e^{i\ell s}}{\sqrt{2\pi}}\Big|_{E_0}\st \ \ell\in
\mathbb{Z}\right\} = \left\{\frac{e^{2^ni\ell
s}}{\sqrt{2\pi}}\Big|_{2^{-n}E_0} \st \ \ell\in
\mathbb{Z}\right\}
\end{equation}
is an orthonormal basis for $L^2(2^{-n}E_0)$ for each $n$. It
follows that
\[
\{D^n\widehat T^\ell \widehat\psi_s\st \ n,\ell\in \mathbb{Z}\}
\]
is an orthonormal basis for $L^2(\mathbb{R})$. Hence
$\{D^nT^\ell\psi_s\st \ n,\ell\in \mathbb{Z}\}$ is an orthonormal
basis for $L^2(\mathbb{R})$, as required.
\subsection{Spectral Set Condition}
From the argument above describing why Shannon's wavelet is, indeed,
an orthonormal basis generator, it is clear that {\em sufficient\/}
conditions for $E$ to be a wavelet set are
\begin{quote}
(i)~~the normalized exponential $\frac1{\sqrt{2\pi}} e^{i\ell s}$, $\ell\in
\mathbb{Z}$, when restricted to $E$ should constitute an orthonormal basis for
$L^2(E)$ (in other words $E$ is a {\em spectral set\/} for the integer lattice
$\mathbb{Z}$),
\end{quote}
\noindent and
\begin{quote}
(ii)~~The family $\{2^nE\st \ n\in\mathbb{Z}\}$ of dilates of $E$ by integral
powers of 2 should constitute a measurable partition (i.e.\ a partition modulo
null sets) of $\mathbb{R}$.
\end{quote}
\noindent These conditions are also necessary. In fact if a set $E$ satisfies
(i), then for it to be a wavelet set it is obvious that (ii) must be satisfied.
To show that (i) must be satisfied by a wavelet set $E$, consider the vectors
\[
\widehat D^n \widehat\psi_E = \frac1{\sqrt{2\pi}} \chi_{2^{-n}E},\qquad n\in
\mathbb{Z}.
\]
Since $\widehat\psi_E$ is a wavelet these must be orthogonal, and so the sets
$\{2^nE\st \ n\in~\mathbb{Z}\}$ must be disjoint modulo null sets. It follows
that $\{\frac1{\sqrt{2\pi}} e^{i\ell s}|_E\st \ \ell\in \mathbb{Z}\}$ is not
only an orthonormal set of vectors in $L^2(E)$, it must also {\em span\/}
$L^2(E)$.
It is known from the theory of \emph{spectral sets} (as an
elementary special case) that a measurable set $E$ satisfies (i) if
and only if it is a generator of a measurable partition of
$\mathbb{R}$ under translation by $2\pi$ (i.e.\ iff $\{E+2\pi
n\st \ n\in \mathbb{Z}\}$ is a measurable partition of
$\mathbb{R}$). This result generalizes to spectral sets for the
integral lattice in $\mathbb{R}^n$. For this elementary special case
a direct proof is not hard.
\subsection{Translation and Dilation Congruence}
We say that measurable sets $E,F$ are {\em translation congruent modulo\/}
$2\pi$ if there is a measurable bijection $\phi\st \  E\to F$ such that
$\phi(s)-s$ is an integral multiple of $2\pi$ for each $s\in E$; or
equivalently, if there is a measurable partition $\{E_n\st \ n\in
\mathbb{Z}\}$ of $E$ such that
\begin{equation}\label{eq39}
\{E_n  + 2n\pi\st \ n\in \mathbb{Z}\}
\end{equation}
is a measurable partition of $F$. Analogously, define measurable sets $G$ and
$H$
to be {\em dilation congruent modulo\/} 2 if there is a measurable bijection
$\tau\st \ G\to H$ such that for each $s\in G$ there is an integer $n$,
depending on $s$, such that $\tau(s) = 2^ns$; or equivalently, if there is a
measurable partition $\{G_n\}^{+\infty}_{n=-\infty}$ of $G$ such that
\begin{equation}\label{eq40}
\{2^nG\}^{+\infty}_{n=-\infty}
\end{equation}
is a measurable partition of $H$. (Translation and dilation congruency modulo
other positive numbers of course make sense as well.)
The following lemma is useful. It is Lemma 4.1 of \cite{DL}.
\begin{lemma}\label{lem7}
Let $f\in L^2(\mathbb{R})$, and let $E = \text{\rm supp}(f)$. Then $f$ has the
property that
\[
\{e^{i\ell s}f\st \ell\in \mathbb{Z}\}
\]
is an orthonormal basis for $L^2(E)$ if and only if
\begin{itemize}
\item[(i)] $E$ is congruent to $[0,2\pi)$ modulo $2\pi$, and
\item[(ii)] $|f(s)| = \frac1{\sqrt{2\pi}}$ a.e.\ on $E$.
\end{itemize}
\end{lemma}
We include a sketch of the proof of Lemma \ref{lem7} for the case  $f = \chi_E$ , because the ideas in it are relevent to the sequel. For general $f$ the proof is a slight modification of this.  If $E$ is a measurable set which is $2\pi$--translation congruent to $[0,2\pi)$, then since
\[
\left\{\frac{e^{i\ell s}}{\sqrt{2\pi}}\Big|_{[0,2\pi)}\Bigg\vert\, n\in\mathbb{Z}\right\}
\]
is an orthonormal basis for $L^2[0,2\pi]$ and the exponentials $e^{i\ell s}$ are $2\pi$--invariant, as in the case of Shannon's wavelet, it follows that
\[
\left\{\frac{e^{i\ell s}}{\sqrt{2\pi}}\Big|_E\Bigg\vert\,  \ell\in \mathbb{Z}\right\}
\]
is an orthonormal basis for $L^2(E)$.
\par
Conversely, if $E$ is not $2\pi$--translation congruent to $[0, 2\pi)$, then either $E$ is congruent to a proper subset  $ \Omega$  of $[0, 2\pi)$, which is not of full measure, or there exists an integer  $k$ such that $E  \cap (E + 2\pi k)$ has positive measure.  In the first case, since the exponentials  $\frac{e^{i\ell s}}{\sqrt{2\pi}}$ restricted to $[0.2\pi)$ do not form an orthonormal basis for $L^2([0, 2\pi))$, the same exponentials restricted to $E$ cannot form an orthonormal basis for $L^2(E)$.  In the second case, let $E_1 = E \cap (E + 2\pi k)$, and $E_2 = E - 2\pi k = (E - 2\pi k) \cap E$, and let $h = \chi_{E_1} - \chi_{E_2}$ .  Then $h \in L^2(E)$ ,  and
$$
h  \perp  \frac{e^{i\ell s}}{\sqrt{2\pi}}\Big|_E,
$$
for all $\ell \in \mathbb{Z}$.  Thus the exponentials $\frac{e^{i\ell s}}{\sqrt{2\pi}}\Big|_E $ cannot even span $L^2(E)$.  This completes the proof sketch.
\par
Next, observe that if $E$ is $2\pi$--translation congruent to $[0,2\pi)$, then since
\[
\{[0,2\pi) + 2\pi n\st \ n\in \mathbb{Z}\}
\]
is a measurable partition of $\mathbb{R}$, so is
\[
\{E + 2\pi n\st \ n\in \mathbb{Z}\}.
\]
Similarly, if $F$ is  $2$--dilation congruent to the Shannon set $E_0 = [-2\pi, -\pi) \cup [\pi, 2\pi)$, then since $\{2^nE_0 : n \in \mathbb{Z}\}$ is a measurable partition of $\mathbb{R}$, so is $\{2^n F : n \in \mathbb{Z}\}$.
\par
These arguments can be reversed.
We say that a measurable subset $G\subset \mathbb{R}$ is a 2--{\em dilation
generator\/} of a {\em partition\/} of $\mathbb{R}$ if the sets
\begin{equation}\label{eq41}
2^nG := \{2^ns\st \ s\in G\},\qquad n\in \mathbb{Z}
\end{equation}
are disjoint and $\mathbb{R}\backslash \cup_n 2^nG$ is a null set. Also, we say
that $E\subset \mathbb{R}$ is a $2\pi$--{\em translation generator of a
partition\/} of $\mathbb{R}$ if the sets
\begin{equation}\label{eq42}
E + 2n\pi := \{s + 2 n\pi\st \ s\in E\},\qquad n\in \mathbb{Z},
\end{equation}
are disjoint and $\mathbb{R}\backslash \cup_n (E+2n\pi)$ is a null set.
\par
The following is Lemma 4.2 of \cite{DL}.
\begin{lemma}\label{lem8}
A measurable set $E\subseteq \mathbb{R}$ is a $2\pi$--translation generator of a
partition of $\mathbb{R}$ if and only if, modulo a null set, $E$ is translation
congruent to $[0,2\pi)$ modulo $2\pi$. Also, a measurable set $G\subseteq
\mathbb{R}$ is a 2-dilation generator of a partition of $\mathbb{R}$ if and only
if, modulo a null set, $G$ is a dilation congruent modulo 2 to the set $[-2\pi,
-\pi) \cup [\pi,2\pi)$.
\end{lemma}
\begin{definition}\label{def999}
By a fundamental domain for a group of (measurable) transformations $\mathcal{G}$ on a measure space $(\Omega,\mu)$ we will mean a measurable set
$C$ with the property that $\{g(C) : g\in\mathcal{G}\}$ is a measurable partition (tessellation) of $\Omega$; that is, $\Omega \setminus \left(\displaystyle{\bigcup_{g\in\mathcal{G}}} g(C)\right)$ is a $\mu$-null set and $g_1 (C) \cap g_2 (C)$ is a $\mu$-null set for $g_1\neq g_2$.
\end{definition}
Thus the sets $E$ and $F$ in Lemma \ref{lem8} are fundamental domains for the groups generated by $2\pi$ and dilation by 2, respectively. Moreover,
Lemma \ref{lem8} actually characterizes the fundamental domains for those groups.
\subsection{A Criterion}
The following is a useful criterion for
wavelet sets. It was published independently by Dai--Larson in \cite{DL}
and by Fang and Wang in \cite{FW} at
about the same time (December 1994).
\begin{proposition}\label{pro9}
Let $E\subseteq\mathbb{R}$ be a measurable set. Then $E$ is a wavelet set if and
only if $E$ is both a 2--dilation generator of a partition (modulo null sets) of
$\mathbb{R}$ and a $2\pi$--translation generator of a partition (modulo null
sets) of $\mathbb{R}$. Equivalently, $E$ is a wavelet set if and only if $E$ is
both translation congruent to $[0,2\pi)$ modulo $2\pi$ and dilation congruent to
$[-2\pi,-\pi) \cup [\pi,2\pi)$ modulo 2. In the terminology of Definition \ref{def999}, $E$ is a wavelet set if and only if $E$ is a fundamental
domain for this dilation group and at the same time a fundamental domain for this translation group.
\end{proposition}
Note that a set is $2\pi$--translation congruent to $[0,2\pi)$ iff it is
$2\pi$--translation congruent to $[-2\pi, \pi)\cup [\pi,2\pi)$. So the last
sentence of Proposition \ref{pro9} can be stated:\ A measurable set $E$ is a
wavelet set if and only if it is both $2\pi$--translation and 2--dilation
congruent to the Littlewood--Paley set $[-2\pi, -\pi)\cup [\pi,2
\pi)$.
\vskip 5pt
For our later purposes, we need the generalization of the above results to $\mathbb{R}^n$. To this end, a few definitions are necessary.
\par
Let $X$ be a metric space and $m$ a $\sigma$-finite non-atomic Borel measure on $X$ for which the measure of every open set is positive and for which
bounded sets have finite measure. Let $\mathcal{T}$ and $\mathcal{D}$ be countable groups of homeomorphisms of $X$ that map bounded sets to bounded
sets and which are absolutely continuously in the sense that they map $m$-null sets to $m$-null sets. Furthermore, let $\mathcal{G}$ be a countable
group of absolutely continuous Borel isomorphisms of $X$. Denote by $\mathcal{B}$ the family of Borel sets of $X$.
\par
The following definition completely generalizes our definitions of $2\pi$--translation congruence and $2$--dilation congruence given in the beginning of subsection 2.4.
\begin{definition}
Let $E, F\in\mathcal{B}$. We call $E$ and $F$ {\em $\mathcal{G}$--congruent} and write $E \sim_{\mathcal{G}} F$, if there exist measurable partitions
$\{E_g : g\in\mathcal{G}\}$ and $\{F_g : g\in\mathcal{G}\}$ of $E$ and $F$, respectively, such that $F_g = g (E_g)$, for all $g\in\mathcal{G}$,
modulo $m$-null sets.
\end{definition}
\begin{proposition}\label{prop99} ${}$
\begin{enumerate}
\item $\mathcal{G}$--congruence is an equivalence relation on the family of $m$-measurable sets.
\item If $E$ is a fundamental domains for $\mathcal{G}$, then $F$ is a fundamental domain for $\mathcal{G}$ iff $F \sim_{\mathcal{G}} E$.
\end{enumerate}
\end{proposition}
\begin{proof}
See \cite{DL}.
\end{proof}
\begin{definition}\label{dil-trans}
We call $(\mathcal{D},\mathcal{T})$ an {\em abstract dilation--translation pair} if
\begin{enumerate}
\item   For each bounded set $E$ and each open set $F$ there exist elements $\delta\in\mathcal{D}$ and $\tau\in\mathcal{T}$ such that $\tau(F)
        \subset\delta(E)$.
\item   There exists a fixed point $\theta\in X$ for $\mathcal{D}$ with the property that if $N$ is any neighborhood of $\theta$ and $E$ any bounded
        set, there is an element $\delta\in\mathcal{D}$ such that $\delta(E)\subset N$.
\end{enumerate}
\end{definition}
The following result and its proof can be found in \cite{DLS1}.
\begin{theorem}\label{t1}
Let $X$, $\mathcal{B}$, $m$, $\mathcal{D}$, and $\mathcal{T}$ as above. Let $(\mathcal{D},\mathcal{T})$ be an abstract dilation--translation pair with $\theta$ being the $\mathcal{D}$ fixed point. Assume that $E$ and $F$ are bounded measurable sets in $X$ such that $E$ contains a neighborhood of $\theta$, and $F$ has non-empty interior and is bounded away from $\theta$. Then there exists a measurable set $G\subset X$, contained in
$\displaystyle{\bigcup_{\delta\in\mathcal{D}}} \delta(F)$, which is both $\mathcal{D}$--congruent to $F$ and $\mathcal{T}$--congruent to $E$.
\end{theorem}
The following is a consequence of Proposition \ref{prop99} and Theorem \ref{t1} and is the key to obtaining wavelet sets.
\begin{corollary}
With the terminology of Theorem \ref{t1}, if in addition $F$ is a fundamental domain for $\mathcal{D}$ and $E$ is a fundamental domain for
$\mathcal{T}$, then there exists a set $G$ which is a common fundamental domain for both $\mathcal{D}$ and $\mathcal{T}$.
\end{corollary}
In order to apply the above result to wavelet sets in $\mathbb{R}^n$, we make the following definition.
\begin{definition}
A \em{dilation $A$--wavelet set} is a measurable subset $E$ of $\mathbb{R}^n$ for which the inverse Fourier transform of $(m(E))^{-1/2}\,\chi_E$ is an orthonormal dilation $A$--wavelet.
\end{definition}
Two measurable subsets $H$ and $K$ of $\mathbb{R}^n$ are called {\em $A$--dilation congruent}, in symbols $H\sim_{\delta_A} K$, if there exist
measurable partitions $\{H_\ell \st \ell\in\mathbb{Z}\}$ of $H$ and $\{K_\ell \st \ell\in\mathbb{Z}\}$ of $K$ such that $K_\ell = A^\ell H_\ell$
modulo Lebesgue null sets. Moreover, two measurable sets $E$ and $F$ of $\mathbb{R}^n$ are called {\em $2\pi$--translation congruent}, written $E\sim_{\tau_{2\pi}} F$, if there exists measurable partitions $\{E_\ell \st \ell\in\mathbb{Z}^n\}$ of $E$ and $\{F_\ell \st \ell\in\mathbb{Z}^n\}$ of $F$ such that $F_\ell = E_\ell + 2\pi\ell$ modulo Lebesgue null sets.
\par
Note that this generalizes to $\R^n$ our definition of $2\pi$--translation congruence for subsets of $\R$.  Observe that $A$--dilation by an expansive matrix together with $2\pi$--translation congruence is a special case (in fact, it is really the {\em prototype} case) of an abstract dilation-translation pair (Definition \ref{dil-trans}).  Let $\mathcal{D}$ be the group of dilations by powers of $A$, $\{A^\ell \st \ell \in \Z\}$ on $\R^n$, and let $\mathcal{T}$ be the group of translations by the vectors $\{2\pi k \st k \in \Z^n\}$.  Let $E$ be any bounded set, and let $F$ be any open set that is bounded away from 0. Let $r > 0$ be such that  $E \subseteq B_r(0)$.  Since $A$ is expansive there is an $\ell \in \mathbb{N}$ such that $A^\ell F$ contains a ball $B$ of radius large enough so that $B$ contains some lattice point  $2k\pi$ together with the ball $B_R(2k\pi)$ of radius $R > 0$ centered at the lattice point.  Then $E + 2k\pi \subseteq A^\ell F$.  That is, the $2k\pi$--translate of $E$ is contained in the $A^\ell$--dilate of $F$,  as required in (1) of Definition \ref{dil-trans}. For (2) of Definition \ref{dil-trans}, let $\theta = 0$, and let $N$ be a neighborhood of 0, and let $E$ be any bounded set.  As above, choose $r  > 0$ with  $E \subseteq B_r(0)$. Let $\ell \in \mathbb{N}$ be such that $A^\ell N$ contains $B_r(0)$.  Then $A^{-\ell}$ is the required dilation such that $A^{-\ell}E \subseteq N$.
\par
Note that if $W$ is a measurable subset of $\mathbb{R}^n$ that is $2\pi$--translation congruent to the $n$-cube $E := \displaystyle{\mathsf{X}_
{k=1}^n} [-\pi, \pi)$, it follows from the exponential form of $\widehat{T}_j$ that $\left\{\widehat{T}_1^{\ell_1}\widehat{T}_2^{\ell_2}\cdots
\widehat{T}_n^{\ell_n}\,(m(W))^{-1/2}\,\chi_W\st \ell = (\ell_1,\ell_2,\ldots,\ell_n)\in\mathbb{Z}^n\right\}$ is an orthonormal basis for $L^2 (W)$.
Furthermore, if $A$ is an expansive matrix, i.e., $A$ is similar to a strict dilation, and $B$ the unit ball of $\mathbb{R}^n$ then with $F_A := A(B)
\setminus B$ the collection $\{A^k F_A : k\in\mathbb{Z}\}$ is a partition of $\mathbb{R}^n\setminus\{0\}$. As a consequence, $L^2 (F_A)$, considered
as a subspace of $L^2 (\mathbb{R}^n)$, is a complete wandering subspace for $D_A$. Hence, $L^2 (\mathbb{R}^n)$ is a direct sum decomposition of the
subspaces $\{D_A^k L^2 (F_A) \st k\in\mathbb{Z}\}$. Clearly, any other measurable set $F^\prime\sim_{\delta_A} F_A$ has this same property.
\par
The above theorem not only gives the existence of wavelet sets in $\mathbb{R}^n$, but also shows that there are sufficiently many to generate the
Borel structure of $\mathbb{R}^{n}$. For details, we refer the reader to \cite{DLS1}. For our purposes, we only quote Corollary 1 in \cite{DLS1},
as a theorem.
\begin{theorem}\label{Th2.11}
Let $n\in\mathbb{N}$ and let $A$ be an expansive $n\times n$ matrix. Then there exist dilation--$A$ wavelet sets.
\end{theorem}
Some concrete examples of wavelet sets in the plane were subsequently obtained by Soardi and Weiland, and others were obtained by Gu and by Speegle
in their thesis work at Texas A\&M University. Two additional examples were constructed by Dai for inclusion in the revised concluding remarks
section of \cite{DL}.
\section{Iterated function systems and fractal functions}
Fractal (interpolation) functions were first systematically introduced in \cite{B}. Their construction is based on iterated function systems and
their properties \cite{BD,Hu}.
\par
To this end, recall that a {\em contraction} on a metrizable space $(M,d)$ is a mapping $f:M\to M$ such that
there exists a $0\leq k < 1$, called the {\em contractivity constant}, so that for all $x,y\in M$
\[
d(f(x),f(y)) \leq k\,d(x,y).
\]
\begin{definition}
Let $(X,d)$ be a complete metrizable space with metric $d$ and let $\{T_i \st i = 1, \ldots, N\}$ be a finite set of contractions on $X$. The pair
$((X,d),\{T_i\})$ is called an iterated function system (IFS) on $X$.
\end{definition}
With the finite set of contractions, one can associate a set-valued operator $\mathscr{T}$, called the Hutchinson operator, defined on the hyperspace
${H}(X)$ of nonempty compact subsets of $X$ endowed with the Hausdorff metric $d_H$:
\[
\mathscr{T}(E) := \bigcup_{i=1}^N T_i (E).
\]
It is easy to show that the Hutchinson operator is contractive on the complete metric space $({H}(X), d_H)$ with contractivity constant $\max_{1\leq
i\leq N} s_i$, where $s_i$ is the contractivity constant of $T_i$. By the Banach Fixed Point Theorem, $\mathscr{T}$ has a unique fixed point, called
the \textit{fractal $F$ associated with the IFS $((X,d),\{T_i\})$}. The fractal $F$ satisfies
\be\label{fixedpoint}
F = \mathscr{T}(F) = \bigcup_{i=1}^N T_i (F),
\ee
i.e., $F$ is made up of a finite number of images of itself. The proof of the Banach Fixed Point Theorem shows that the fractal can be iteratively
obtained via the following procedure. Choose $F_0\in {H}(X)$ arbitrary. Define
\[
F_n := \mathscr{T}(F_{n-1}),\qquad n\in\mathbb{N}.
\]
Then $F = \displaystyle{\lim_{n\to\infty}} F_n$, where the limit is taken in the Hausdorff metric.
\par
A special situation occurs when $X := [a,b]\times\mathbb{R}\subset\mathbb{R}^2$, $a < b$, and $d$ is the Euclidean metric. Let $\{(x_j,y_j) \st x_0
:= a < x_1 < \ldots x_N := b,\,y_j\in\mathbb{R},\,j = 0,1,\ldots, N\}$ be a given set of interpolation points. For $i = 1, \ldots, N$, let $s_i\in
(-1,1)$ and let $F_0:= [a,b]\times [a,b]$. Define images $T_i F_0$, $i = 1, \ldots, N$, of $F_0$ as follows. $T_i F_0$ is the unique parallelogram
with vertices at $(x_{i-1},y_{i-1})$, $(x_i,y_i)$, $(x_i,y_i + s_i (b-a))$, and $(x_{i-1},y_{i-1} + s_i (b-a))$. There exists a unique affine mapping
$T_i : X\to X$ such that $T_i F_0 = T_i (F_0)$, namely
\[
T_i \begin{pmatrix} x\\ y \end{pmatrix} = \begin{pmatrix} a_i & 0\\ c_i & s_i\end{pmatrix}\begin{pmatrix} x\\ y \end{pmatrix} + \begin{pmatrix}
\alpha_i\\ \beta_i \end{pmatrix},
\]
where
\begin{align*}
a_i & := \frac{x_i-x_{i-1}}{b-a}, & c_i &:= \frac{y_i-y_{i-1} - s_i\,(y_N - y_0)}{b-a},\\
\alpha_i & := \frac{b x_{i-1} - a x_i}{b - a}, & \beta_i &:= \frac{b y_{i-1} - a y_i - s_i\,(b y_0 - a y_N)}{b - a}.
\end{align*}
Since the scaling factors $s_i$ are in modulus less than one, the affine mappings $T_i$ are contractive on $X$ and thus $((X,d),\{T_i\})$ is an IFS.
As such, it has a unique fixed point $F$, which turns out to be the graph of a continuous function $f: [a,b]\to \mathbb{R}$ satisfying $f(x_j) =
y_j$, for $j = 0,1,\ldots, N$. (See \cite{B}) The graph of $f$ is in general a fractal set in the above sense and contains the given set of
interpolation points.
\begin{example}
Let $a = 0$, $b =1$, $N = 2$, and choose interpolation points $\{(0,0),(0.5,0.7),(1,0)\}$ and scaling factors $s_1 = 0.6$ and $s_2 = 0.4$. The
sequence of graphs in Figure~\ref{fig1} shows the geometric construction of a fractal function with these parameters as outlined above.
\begin{figure}[h]
{\begin{center}\includegraphics[width=2in,height=1in]{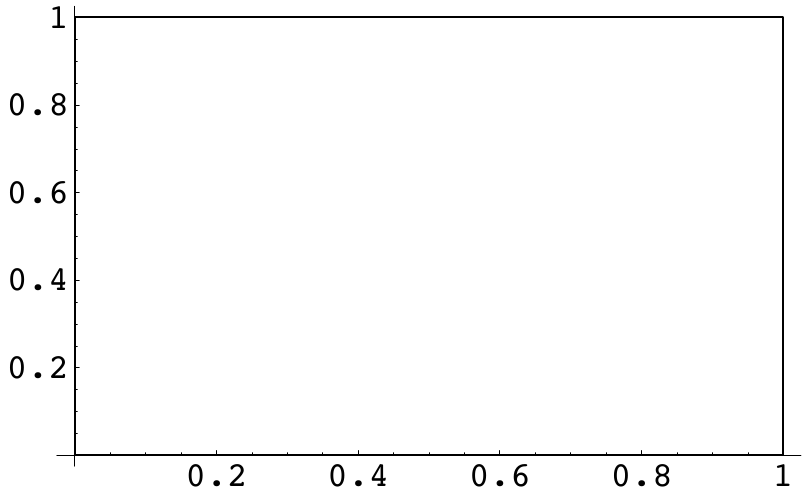}\hspace{1cm}\includegraphics[width=2in,height=1in]{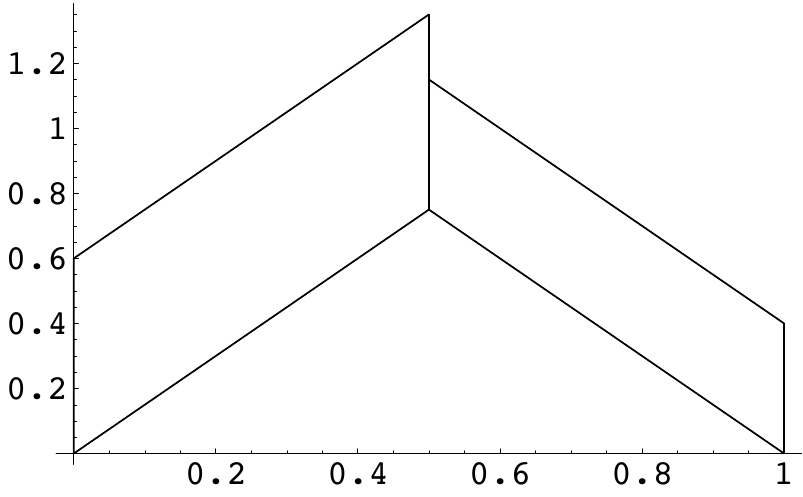}\end{center}}%
{\begin{center}\includegraphics[width=2in,height=1in]{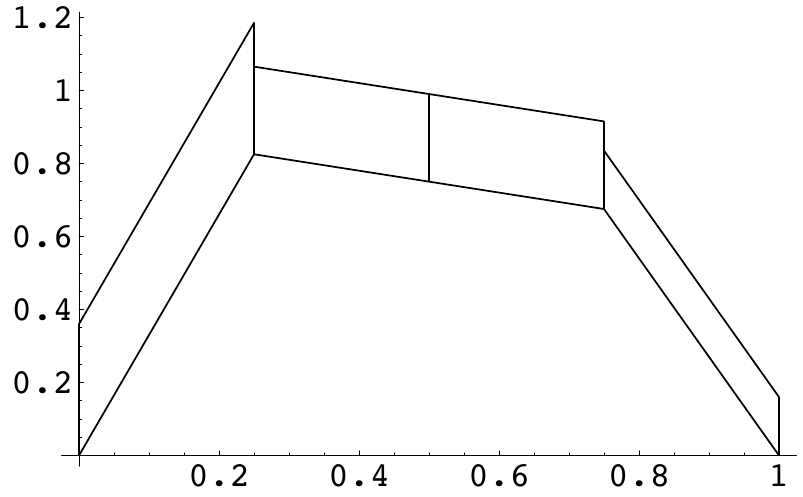}\hspace{1cm}\includegraphics[width=2in,height=1in]{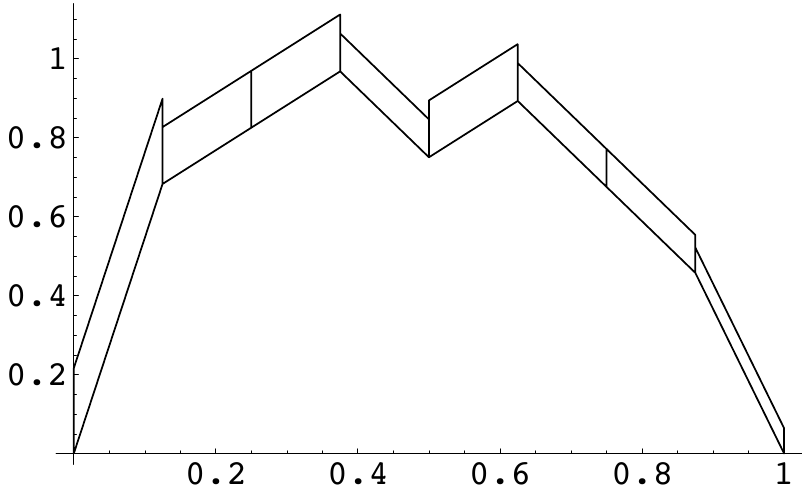}\end{center}}%
{\begin{center}\includegraphics[width=2in,height=1in]{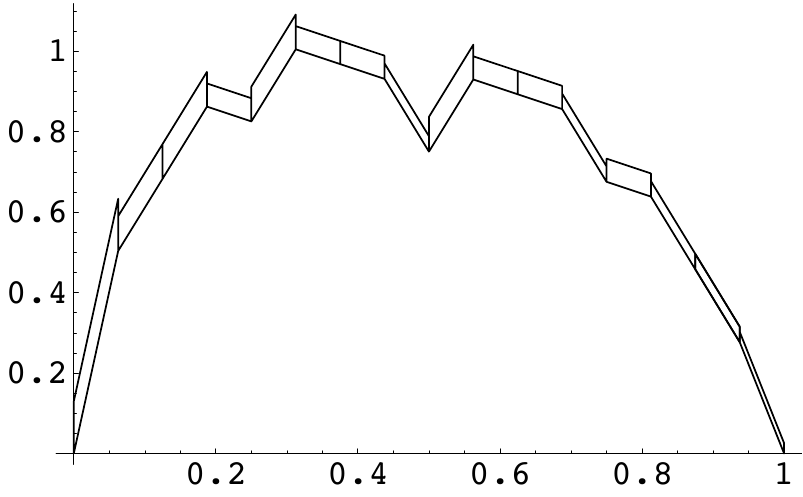}\hspace{1cm}\includegraphics[width=2in,height=1in]{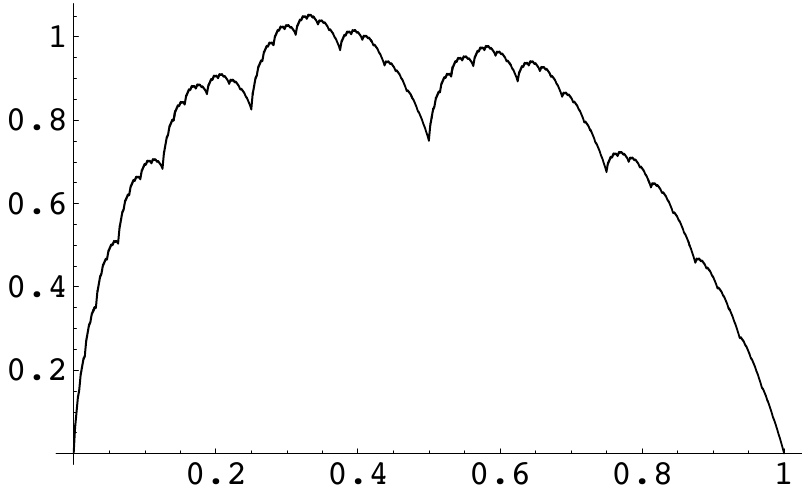}\end{center}}%
\caption{The geometric construction of a fractal function}\label{fig1}
\end{figure}
\end{example}
Writing the fixed point equation~(\ref{fixedpoint}) using the affine mappings $T_i$ for a point $(x, f(x))$ on the graph of a continuous fractal
function $f$ yields
\[
\begin{pmatrix} x\\ f(x) \end{pmatrix}\Bigg\vert_{x\in [x_{i-1},x_i]} = \begin{pmatrix} a_i & 0\\ c_i & s_i\end{pmatrix}
\begin{pmatrix} x\\ f(x)\end{pmatrix}\Bigg\vert_{x\in [a,b]} + \begin{pmatrix} \alpha_i\\ \beta_i \end{pmatrix},
\]
Setting $u_i (x) := a_i x + \alpha_i$ and $p_i (x) := c_i x + \beta_i$, $i = 1,\ldots, N$, one can rewrite the second component of the above equation
as
\[
f(x) = p_i(u_i^{-1}(x)) + s_i f(u_i^{-1}(x)), \qquad x\in [x_{i-1},x_i],
\]
or, equivalently,
\be\label{RB}
\begin{split}
f(x) &= \sum_{i=1}^N \left[p_i (u_i^{-1}(x)) + s_i f(u_i^{-1}(x))\right]\chi_{[x_{i-1},x_i]}\\
     & =: p(x) + \sum_{i=1}^N s_i f(u_i^{-1}(x))\chi_{[x_{i-1},x_i]},\qquad x\in [a,b],
\end{split}
\ee
where $\chi$ is the indicator function and $p$ the linear spline whose restriction to any of the intervals $[x_{i-1},x_i]$ equals $p_i\circ u_i^{-1}$. Such a fractal function may be considered as a {\em linear spline parametrized by the row vector $\boldsymbol{s} := (s_1,\ldots, s_n)$.}
\par
Note that the functions $p_i$ are uniquely determined by the interpolation points and that the fractal functions is uniquely determined by the row
vector $\boldsymbol{p}:=(p_1,\ldots, p_N)$ and the scaling factors $(s_1,\ldots, s_n)$. We suppress the dependence of a fractal function $f$ on
$(s_1,\ldots, s_n)$ but write $f = f_{\boldsymbol{p}}$ when necessary.
\par
Denote by $\Pi^1 [a,b]$ the linear space of real polynomials of degree at most one on $[a,b]$ and by $C[a,b]$ the space of continuous functions on
$[a,b]$.
\begin{theorem}
The mapping $\Theta: \Pi^1[a,b]^N\cap C[a,b]\ni\boldsymbol{p}\mapsto f_{\boldsymbol{p}}$ is a linear isomorphism.
\end{theorem}
\begin{proof}
The result follows from the above observations and the uniqueness of the fixed point: $f_{\alpha\boldsymbol{p} + \boldsymbol{q}} = \alpha\,
f_{\boldsymbol{p}} + f_{\boldsymbol{q}}$, $\alpha\in\mathbb{R}$ and $\boldsymbol{p},\boldsymbol{q}\in \Pi^1[a,b]^N$. (See also \cite{M2,M3}.)
\end{proof}
Since the space $\Pi^1[a,b]^N\cap C[a,b]$ is $(N+1)$-dimensional ($2N$ free parameters plus $N-1$ join-up conditions at the interior interpolation
points), the space $\mathfrak{F}^1[a,b]$ of all continuous fractal functions on $[a,b]$ generated by functions in $\Pi^1[a,b]^N$ has also dimension
$N+1$. An basis for $\mathfrak{F}^1[a,b]$ can be found by choosing as the $i$th basis function the continuous fractal function $e_i$ that
interpolates according to
\[
e_i (x_j) = \delta_{ij},\quad i,j = 0,1,\ldots, N.
\]
It follows immediately from the uniqueness of the fixed point, that every continuous fractal function $f\in\mathfrak{F}[a,b]$ can be written as a
linear combination of the form
\be
f(x) = \sum_{j=0}^N y_j \, e_j (x).
\ee
Using the Gram-Schmidt Orthonormalization procedure and the fact that the $L^2$-inner product of two fractal functions in $\mathfrak{F}^1[a,b]$ over the same knot set $\{x_j\st j = 0,1,\dots, n\}$ can be explicitly computed in terms of the parameters $(s_1, \ldots, s_N)$ and the functions in
$\Pi^1[a,b]^N\cap C[a,b]$ (cf. \cite{M2}), the basis $\{e_j\}$ can also be orthonormalized.
\par
Equation~(\ref{RB}) can be interpreted as the fixed point equation for a function-valued operator $\mathscr{B}$ and this interpretation leads to a
more abstract definition of fractal functions. The following theorem gives a general construction of fractal functions in terms of so-called
Read-Bajraktarevi\'{c} operators.
\begin{theorem}
Let $\Omega\subset\mathbb{R}$ be compact and $1< N\in\mathbb{N}$. Assume that $u_i:\Omega\to\Omega$ are contractive homeomorphisms inducing a
partition on $\Omega$, $\lambda_i:\mathbb{R}\to\mathbb{R}$ are bounded functions and $s_i$ real numbers, $i = 1,\ldots, N$.
Let
\be
\mathscr{B} (f) := \sum_{i=1}^{N} \left[\lambda_i\circ u_i^{-1} + s_i\,f\circ u_i^{-1}\right]\chi_{u_i(\Omega)}
\ee
If $\max\{|s_i|\} < 1$, then the operator $\mathscr{B}$ is contractive on $L^\infty (\Omega)$ and its unique fixed point $f: \Omega\to\mathbb{R}$
satisfies
\[
f = \sum_{i=1}^{N} \left[\lambda_i\circ u_i^{-1} + s_i\,f\circ u_i^{-1}\right]\chi_{u_i(\Omega)}
\]
\end{theorem}
\begin{proof}
Apply the Banach Fixed Point Theorem.
\end{proof}
The fixed point of such an operator is called a {\em fractal function}. Note again that $f$ depends on the row  vector of functions $\boldsymbol{\lambda} := (\lambda_1,\ldots, \lambda_N)$.
\par
If $\mathscr{B}$ acts on a normeable or metrizable function space $\mathcal F$, then its fixed point, under appropriate conditions on
$\boldsymbol{\lambda}$ and the scaling factors $(s_1,\ldots, s_n)$, is also an element of $\mathcal F$. In this manner, one can construct fractal
functions with prescribed regularity or approximation properties \cite{M4}.
\par
\begin{figure}[h]
{\begin{center}\includegraphics[width=2in,height=1in]{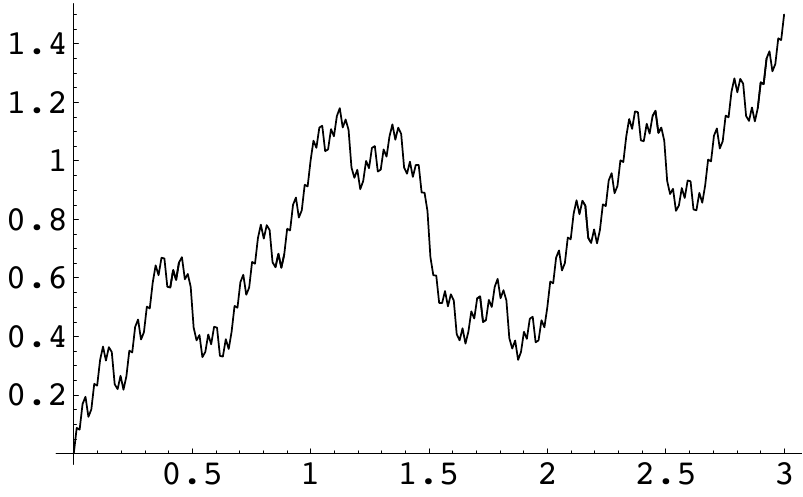}\hspace{1cm}\includegraphics[width=2in,height=1in]{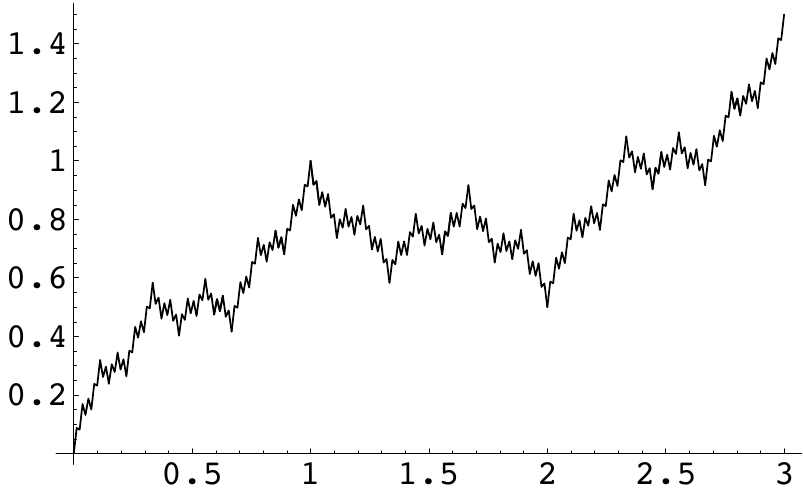}\end{center}}
\caption{A fractal function generated via translations (left) and reflections (right).}\label{fig2}
\end{figure}

In case the contractive homeomorphisms $u_i$, $i=1,\ldots, N$, induce a uniform partition of $\Omega$, the above expressions become more transparent.
In addition, w.l.o.g., suppose that $\Omega = [0,N]$. A natural uniform partition in this case is $[0,N] = \displaystyle{\bigcup_{i=1}^{N-1}} [i-1,i)\,\cup [N-1,N]$ and the mappings $u_i$ are given by
\[
u_i (x) = \frac{x}{N} + i - 1, \qquad i = 1, \ldots, N.
\]
Hence, it suffices to define $u_1 := x/N$ and then all other mappings are given by translating $u_1$: $u_i = u_1 +(i-1) = u_{i-1} + 1$, $i = 2,
\ldots, N$. The fixed point equation for a continuous fractal function than reads
\[
f(x) = \sum_{i=1}^N\left[\lambda_i(N(x - i +1)) + s_i\,f(N(x - i + 1))\right]\chi_{[i-1,i]},\quad x\in [0,N].
\]
Instead of using translations to obtain $u_2,\ldots, u_N$ from $u_1$, one may choose reflections about the partition points $i$ of $[0,N]$, $i = 1,
\ldots, N-1$, instead. The reflection $R_i$ about the point $(i,0)$ on the $x$-axis is given by $R_i(x) = 2i - x$. As above, let $u_1 (x) = x/N$ and
set $u_i = R_{i-1}\circ u_{i-1}$, $i = 2,\ldots, N$. This also generates a uniform partition of $[0,N]$.
\begin{example}
As an example of the two types of continuous fractal functions, we consider $N := 3$ and $s_1 = s_2 = s_3 := 0.5$. The continuous fractal function
$f$ is generated by translations in the above sense with $\boldsymbol{\lambda} = ((1/3 - s_1/2)x, (-1/6 - s_2/2)x + 1, (1/3 - s_3/2)x + 1/2)$,
whereas the continuous fractal function $g$ is generated using reflections and with $\boldsymbol{\lambda} = ((1/3 - s_1/2)x, (1/6 - s_2/2)x + 1, (
1/3 - s_3/2)x + 1/2)$. The values of both functions are the partition points $(i,0)$, $i = 1,\ldots, 4$, are $0,1,1/2$, and $3/2$. The graphs of
these two fractal functions are displayed in Figure~\ref{fig2}.
\end{example}
\par
\begin{figure}[h]
{\begin{center}\includegraphics[width=2in,height=1in]{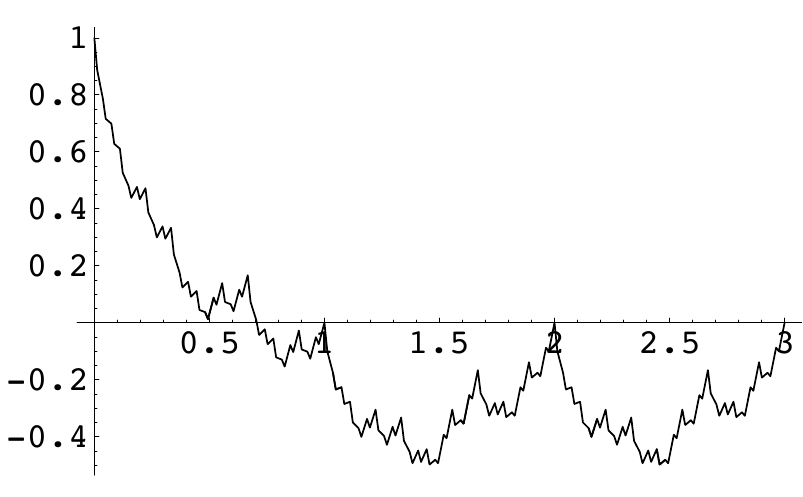}\hspace{1cm}\includegraphics[width=2in,height=1in]{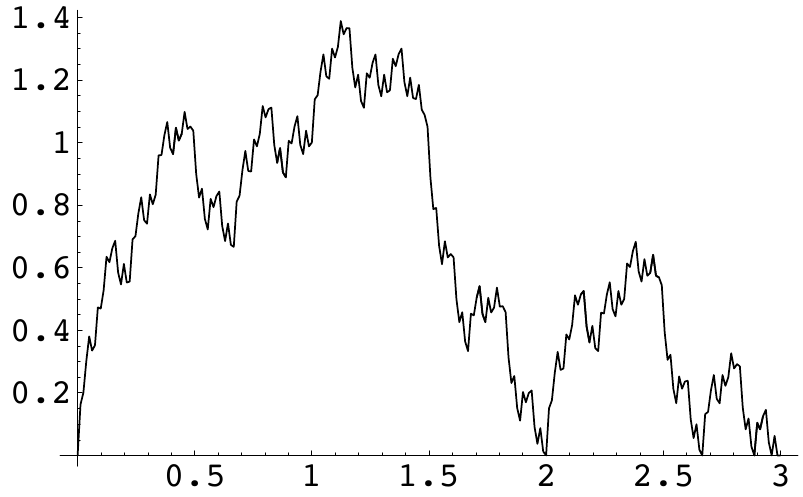}\end{center}}%
{\begin{center}\includegraphics[width=2in,height=1in]{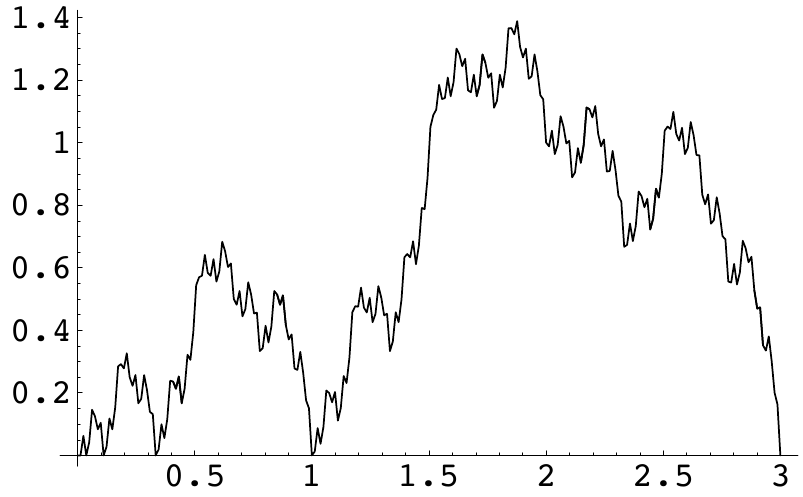}\hspace{1cm}\includegraphics[width=2in,height=1in]{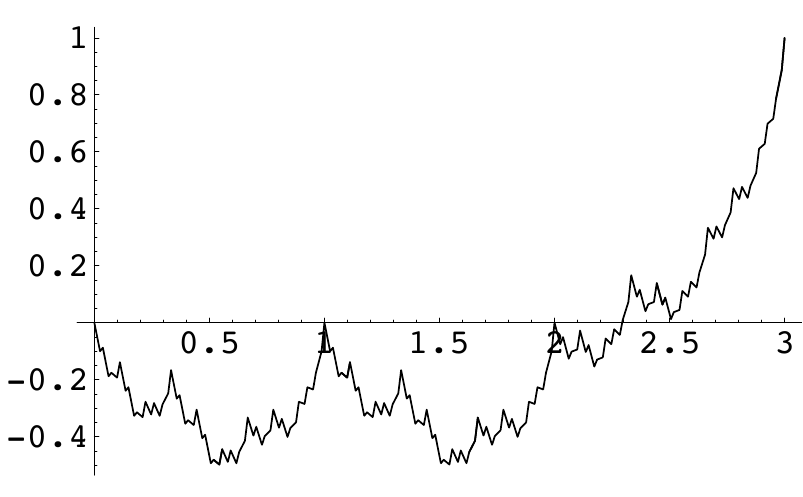}\end{center}}
\caption{The four basis functions for the fractal function $f$.}\label{fig3}
\end{figure}
\begin{figure}[h]
{\begin{center}\includegraphics[width=2in,height=1in]{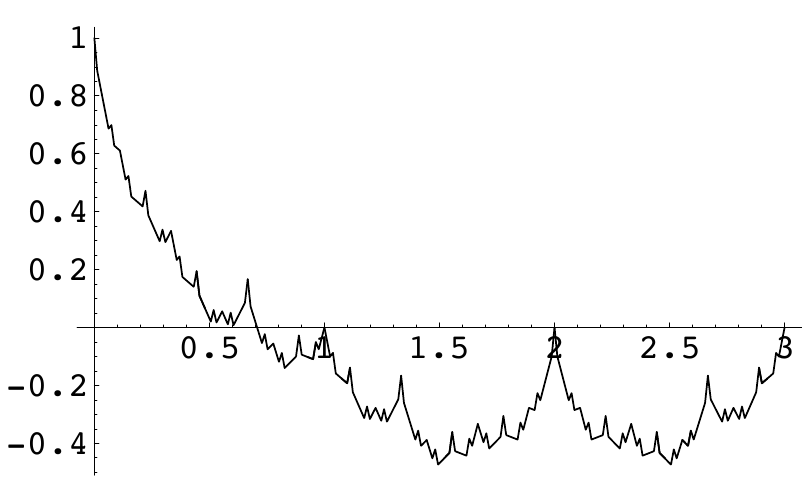}\hspace{1cm}\includegraphics[width=2in,height=1in]{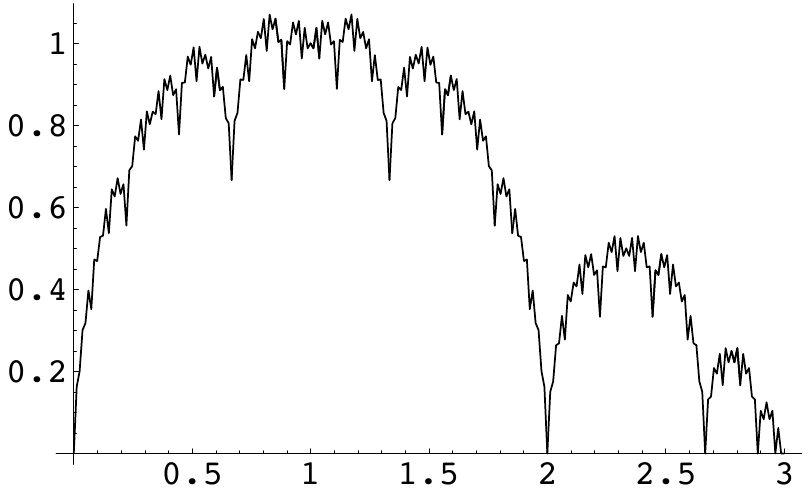}\end{center}}%
{\begin{center}\includegraphics[width=2in,height=1in]{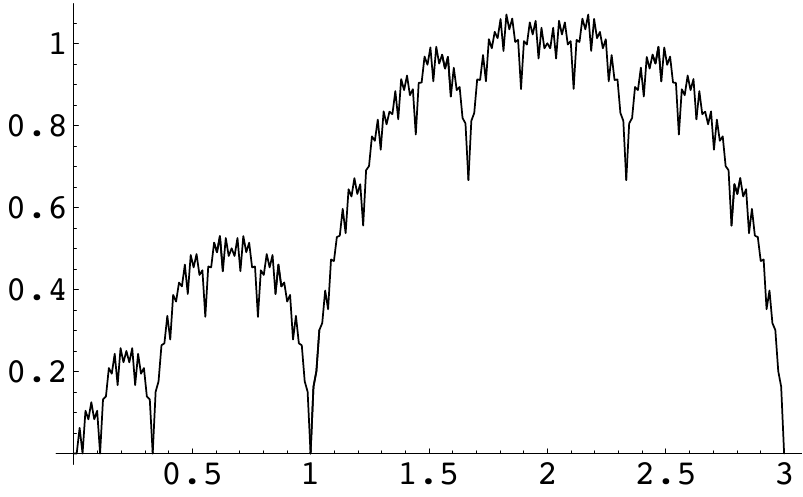}\hspace{1cm}\includegraphics[width=2in,height=1in]{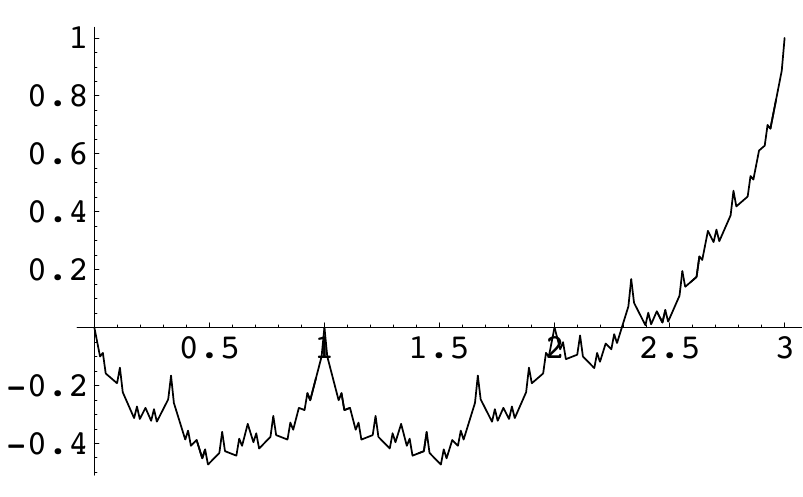}\end{center}}
\caption{The four basis functions for the fractal function $g$.}\label{fig4}
\end{figure}

Both fractal functions belong to the four-dimensional linear space $\Pi^1[a,b]^3\cap C[a,b]$ and the four basis functions are depicted in
Figures~\ref{fig3} and~\ref{fig4}.
\par
At this point the question arises whether there is a generalization of the above procedures to higher dimensions. It is possible to define fractal
surfaces in $\mathbb{R}^n$, $1 < n\in\mathbb{N}$, using reflections instead of translations in the definition of the contractive homeomorphisms
$u_i$. We will see that this leads to a natural way of tessellating the embedding space and to the construction of a multiresolution analysis of
$L^2(\mathbb{R}^n)$.
\section{Coxeter groups and foldable figures}
In order to carry out the construction of fractal surfaces using reflections, a short excursion into the theory of Coxeter groups and foldable
figures is necessary. The interested reader is referred to \cite{Bo,C,G,Gu,H,HW} for more details and proofs.
\subsection{Coxeter groups}
\begin{definition}
A Coxeter group $\mathcal{C}$ is a discrete groups with a finite number of generators $\{r_i\,\st\,i = 1,\ldots, k\}$ satisfying
\[
\mathcal{C} = \bigl\langle r_1,\ldots, r_k\,\st\, (r_i r_j)^{m_{ij}} = 1,\; 1 \leq i,j \leq k\bigr\rangle
\]
where $m_{ii} = 1$, for all $i$, and $m_{ij}\geq 2$, for all $i\neq j$. ($m_{ij}=\infty$ is used to indicate that no relation exists.)
\end{definition}
A geometric representation of a Coxeter group is given by considering it as a subgroup of $GL (V)$, where $V$ is a $k$-dimensional real vector space, which we take to be $\mathbb{R}^k$ endowed with its usual positive definite symmetric bilinear form $\inn{\cdot}{\cdot}$. In this representation, the generators are interpreted in the following way.
\par
A reflection about a linear hyperplane $H$ is defined as a linear mapping $\rho: V \to V$ such that $\rho\vert_{H
} = \textrm{id}_H$ and $\rho(x) = -x$, if $x\in H^\perp$. In other words, $\rho$ is an isometric isomorphism of $V$.
\par
Now suppose that $0\neq r\in H^\perp$, then an easy computation shows that
\[
\rho_r(x) = x - \displaystyle{\frac{2\inn{x}{r}}{\inn{r}{r}}}\,r
\]
is the reflection about the hyperplane $H$ perpendicular to $r$. One can show that the linear mappings $\rho_{r_i}$, where $\{r_i \st i = 1,\ldots,
k\}$ are the generators of a Coxeter group $\mathcal {C}$, satisfy $(\rho_{r_i}\rho_{r_j})^{m_{ij}} = \textrm{id}_V$. It is known that the map
$r_i\mapsto\rho_{r_i}$ extends to a faithful representation of $\mathcal{C}$ into $GL(V)$.
\par
If one considers the group generated by real reflections about linear hyperplanes in $\mathbb{R}^k$, then this group is isomorphic to a finite
Coxeter group whose $k$ generators correspond to the (unit) normal vectors of the set of hyperplanes.
\subsection{Roots systems and Weyl groups}
The normal vectors to a set of hyperplanes play an important role in the representation theory for Coxeter groups. We have seen above that they
correspond to the generators of such groups. Two such normal vectors, $\pm\, r$, that are orthogonal to a hyperplane are called {\em roots}.
\begin{definition}
A \textit{root system} ${\mathcal{R}}$ is a finite set of nonzero vectors $r_1, \ldots, r_k\in$ $\mathbb{R}^n$ satisfying
\begin{enumerate}
\item   $\mathbb{R}^n = \textrm{span}\,\{r_1,\ldots, r_k\}$
\item   $r, \alpha r\in\mathcal{R}$ iff $\alpha = \pm 1$
\item   $\forall r,s\in\mathcal{R}$: $s - \displaystyle{\frac{2\inn{s}{r}}{\inn{r}{r}}}\,r\in\mathcal{R}$, i.e., the root system $\mathcal{R}$ is
        closed with respect to the reflection through the hyperplane orthogonal to $r$.
\item   $\forall r,s\in\mathcal{R}$: $\displaystyle{\frac{2\inn{s}{r}}{\inn{r}{r}}}\in\mathbb{Z}$, i.e., $\rho_r (s) - s \in\mathbb{Z}$
\end{enumerate}
A subset $\mathcal{R}^+\subset\mathcal{R}$ is called a set of \textit{positive} roots if there exists a vector $v\in\mathbb{R}^n$ such that
$\inn{r}{v} > 0$ if $r\in\mathcal{R}^+$, and $\inn{r}{v} > 0$ if $r\in\mathcal{R}\setminus\mathcal{R}^+$. Roots that are not positive are called
\textit{negative}. Since $r$ is negative iff $-r$ is positive, there are exactly as many positive as there are negative roots.
\par
The group generated by the set of reflections $\{\rho_r\,\st\,r\in\mathcal{R}\}$ is called the \textit{Weyl Group} $\mathcal{W}$ of $\mathcal{R}$.
\end{definition}
It follows from the definition of root system, that the Weyl group $\mathcal{W}$ has finite order, indeed it is a finite Coxeter group.
\begin{example}
A simple example of a Weyl group in $\mathbb{R}^2$ is given by the root system depicted in Figure \ref{klein}. The roots are $r_1 = -r_3 =
(1,0)^\top$ and $r_2 = - r_4 = (0,1)^\top$. The positive roots are $r_1$ and $r_2$. The group of reflections generated by these four roots
is given by
\[
V_4 := \bigl\langle \rho_1,\rho_2\,\st\, \rho_1^2 = \rho_2^2 = 1,\; (\rho_1\rho_2)^2 =1\bigr\rangle,
\]
where $\rho_1$ and $\rho_2$ denotes the reflection about the $y$-, respectively, $x$-axis. This group is commutative and called {\em Klein's four-group} or the {\em  group of order four}. In the classification scheme of Weyl groups $V_4$ is referred to as $A_1\times A_1$ since it is the direct product of the group $A_1 := \bigl\langle \rho_1\,\st\, \rho_1^2 = 1\bigr\rangle$ whose root system is $\mathcal{R} = \{r_1, r_3\}$ with itself.
\begin{figure}[h]
\begin{center}
\includegraphics[width=1in,height=1in]{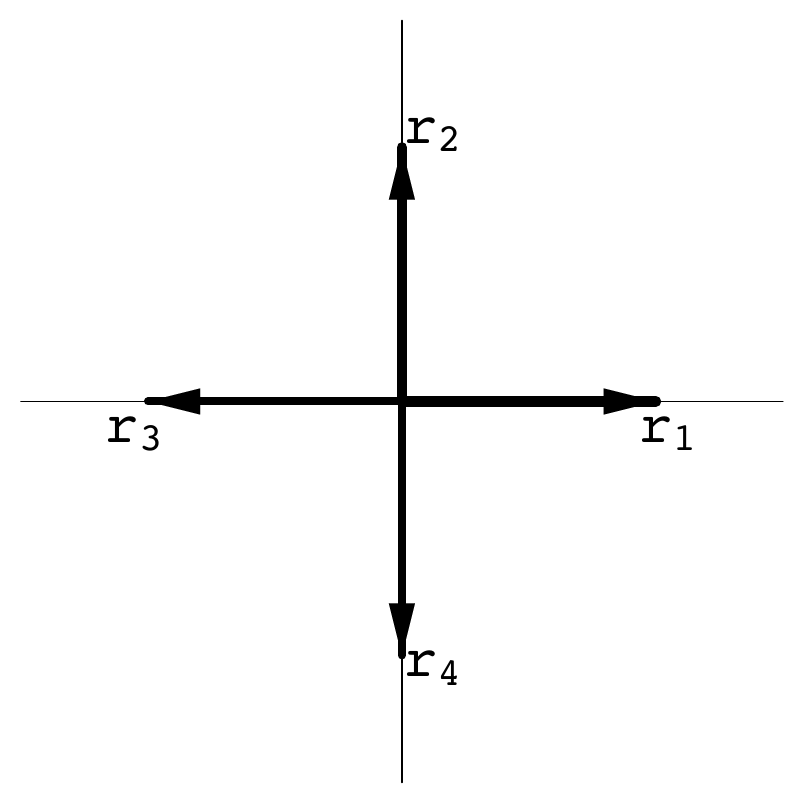}
\end{center}
\caption{The root system for Klein's four-group.}\label{klein}
\end{figure}
\end{example}
\par
For the following, we need some properties of roots systems and Weyl groups, which we state in a theorem.
\begin{theorem}
Let $\mathcal{R}$ be a root system and $\mathcal{W}$ the associated Weyl group. Then the following hold.
\begin{enumerate}
\item   Every root system $\mathcal{R}$ has a basis $\mathcal{B} = \{b_i\}$ consisting of positive (negative) roots.
\item   Let $C_i := \{x\in\mathbb{R}^n\,\st\,\inn{x}{b_i} > 0\}$ be the {\em Weyl chamber} corresponding to the basis $\mathcal{B}$. Then
        the Weyl group $\mathcal{W}$ acts simply transitively on the Weyl chambers.
\item   The set $C := \overline{\bigcap_i C_i}$ is a noncompact fundamental domain for the Weyl group $\mathcal{W}$. It is a simplicial cone, hence
        convex and connected.
\end{enumerate}
\end{theorem}
In order to introduce foldable figures below, we need to consider reflections about affine hyperplanes. For this purpose, let $\mathcal{R}$ be a root
system. An {\em affine hyperplane} with respect to $\mathcal{R}$ is given by
\begin{equation}\label{hyper}
H_{r,k} := \{x\in\mathbb{R}^n\st\inn{x}{r} = k\},\qquad k\in\mathbb{Z}.
\end{equation}
It is easy to show that reflections about affine hyperplanes have the form
\begin{equation}\label{affref}
\rho_{r,k}(x) = x - \displaystyle{\frac{2(\inn{x}{r}-k)}{\inn{r}{r}}}\,r =: \rho_r (x) + k\,{r}^\vee,
\end{equation}
where $r^\vee := 2\,r/\inn{r}{r}$ is the {\em coroot} of $r$.
\begin{definition}
The {\em affine Weyl group} $\widetilde{\mathcal{W}}$ for a root system $\mathcal{R}$ is the (infinite) group generated by the reflections $\rho_{r,k}$ about the affine hyperplanes $H_{r,k}$:
\[
\widetilde{\mathcal{W}} := \bigl\langle \rho_{r,k}\st r\in\mathcal{R}, k\in\mathbb{Z}\bigr\rangle
\]
We sometimes  will refer to the concatenation of elements from $\widetilde{\mathcal{W}}$ as {\em words}.
\end{definition}
\begin{theorem}\label{th4.5}
The affine Weyl group $\widetilde{\mathcal{W}}$ of a root system $\mathcal{R}$ is the semi-direct product $\mathcal{W}\ltimes \Gamma$, where $\Gamma$
is the abelian group generated by the coroots ${r}^\vee$. Moreover, $\Gamma$ is the subgroup of translations of $\widetilde{\mathcal{W}}$ and
$\mathcal{W}$ the isotropy group (stabilizer) of the origin. The group $\mathcal{W}$ is finite and $\Gamma$ infinite.
\end{theorem}
\begin{remark}
There exists a complete classification of all irreducible affine Weyl groups and their associated fundamental domains. These groups are given as
types $A_n$ ($n\geq 1$), $B_n$ ($n\geq 2$), $C_n$ ($n\geq 3$), and $D_n$, ($n\geq 4$), as well as $E_n$, $n = 6, 7, 8$, $F_4$, and $G_2$. (For more
details, we refer the reader to \cite{Bo} or \cite{H}.)
\end{remark}
We need a few more definitions and related results. By a \textit{reflection group} we mean a group of transformations generated by the reflections
about a finite family of affine hyperplanes. Coxeter groups and affine Weyl groups are examples of reflections groups.
\par
Let $\mathcal{G}$ be a reflection group and $\mathcal{O}_n$ the group of linear isometries of $\mathbb{R}^n$. Then there exists a homomorphism
$\phi: \mathcal{G}\to\mathcal{O}_n$ given by
\[
\phi(g)(x) = g(x) - g(0),\quad g\in\mathcal{G},\; x\in\mathbb{R}^n.
\]
The group $\mathcal{G}$ is called essential if $\phi(\mathcal{G})$ only fixes $0\in\mathbb{R}^n$. The elements of $\ker\phi$ are called translations.
\subsection{Foldable figures}
\begin{figure}[h]
\begin{center}
\includegraphics[width=1in,height=1in]{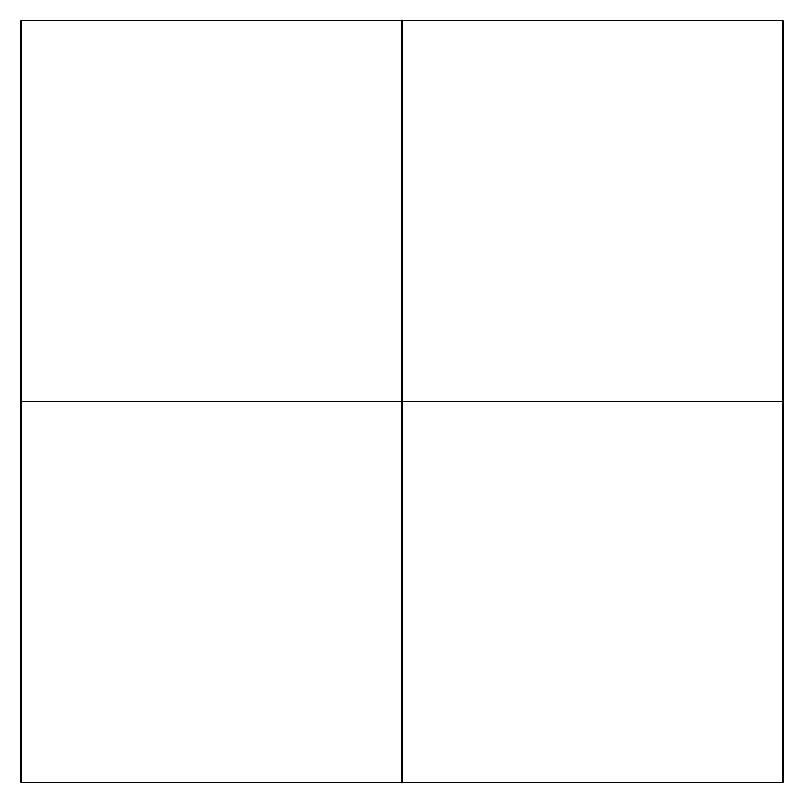}\hspace{2cm}
\includegraphics[width=1in,height=1in]{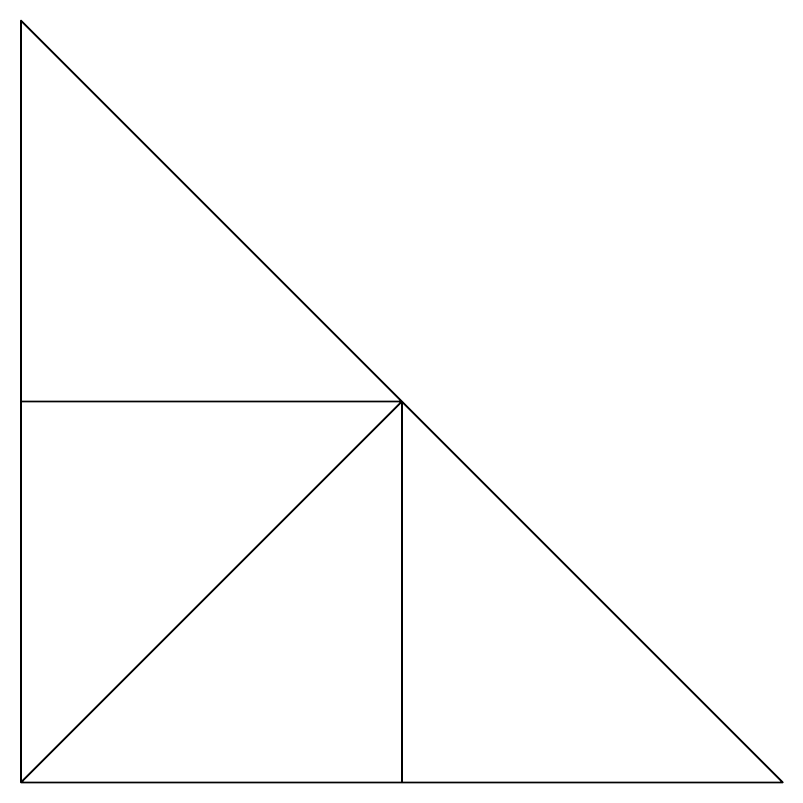}%
\end{center}
\caption{Examples of foldable figures.}\label{fig5}
\end{figure}

In this subsection, we define for our later purposes the important concept of a foldable figure \cite{HW}.
\begin{definition}
A compact connected subset $F$ of $\mathbb{R}^n$ is called a \textit{foldable figure} iff there exists a finite set $\mathcal{S}$ of affine
hyperplanes that cuts $F$ into finitely many congruent subfigures $F_1, \ldots, F_m$, each similar to $F$, so that reflection in any of the cutting
hyperplanes in $\mathcal{S}$ bounding $F_k$ takes it into some $F_\ell$.
\end{definition}
In Figure~\ref{fig5} are two examples of foldable figures shown.
Properties of foldable figures are summarized in the theorem below. The statements and their proofs can be found in \cite{Bo} and \cite{HW}.
\begin{theorem}
${}$
\begin{enumerate}
\item   The reflection group generated by the reflections about the bounding hyperplanes of a foldable figure $F$ is the affine Weyl group
        $\widetilde{W}$ of some root system. Moreover, $\widetilde{W}$ has $F$ as a fundamental domain.
\item   Let $\mathcal{G}$ be a reflection group that is essential and without fixed points. Then $\mathcal{G}$ has a compact fundamental domain.
\item   There exists a one-to-one correspondence between foldable figures and reflection groups that are essential and without fixed points.
\end{enumerate}
\end{theorem}
\section{Fractal surfaces on foldable figures}
Affine fractal surfaces were first systematically introduced in \cite{M1} and slightly generalized in \cite{GH}. A further generalization was
presented in \cite{HM}. The construction is again based on IFS's now defined on simplicial regions $\Delta\subset\mathbb{R}^n$ such that the
contractive homeomorphisms $u_i$ are affine mappings from $\Delta$ into itself. In the present setting, we assume that the domain $\Delta$ is a
foldable figure in the sense of the previous section and therefore partitioned into $N$ congruent subsimplices $\Delta_i$, $i = 1, \ldots, N$:
\[
\Delta = \bigcup_{i=1}^N \Delta_i,\quad\text{and}\quad\overset{\circ}{\Delta}_i\cap\overset{\circ}{\Delta}_j = \varnothing,\;i\neq j,
\]
Thus, one can find $N$ similitudes $u_i: \Delta\to \Delta_i$ of the form
\[
u_i = \sigma\,O_i + b_i
\]
where $O_i\in E(n)$, the Euclidean group of $\mathbb{R}^n$, $b_i\in\mathbb{R}^n$, $i = 1, \ldots, N$, and $\sigma\in (0,1)$ is the similarity ratio
between $\Delta_i$ and $\Delta$. As in Section 2, let $s\in (-1,1)$ be an arbitrary scaling factor and $\{\lambda_i:\mathbb{R}^n\to\mathbb{R} \st\,
i = 1, \ldots, N\}$ a finite collection of continuous affine functions satisfying the following condition.
\[
(*)\quad \begin{cases} \textrm{Let $e_{ij}$ be the common face of $u_i(\Delta)$ and $u_j(\Delta)$. Then $\lambda_i (x) = \lambda_j (x)$}\\
\textrm{for all $x\in u_i^{-1} (e_{ij}) = u_j^{-1}(e_{ij})$, $i,j = 1, \ldots,N$.}\end{cases}
\]
Define an operator $\mathscr{B}: C(\Delta)\to L^\infty (\Delta)$ by
\be\label{10}
\mathscr{B} f (x) := \sum_{i=1}^N \left[\lambda_i \circ u_i^{-1} (x) + s\,f\circ u_i^{-1} (x)\right]\,\chi_{\Delta_i}.
\ee
It can be shown \cite{HM,GH,M1,M2} that $\mathscr{B}$ maps $C(\Delta)$ into itself and is a contractive operator in the $\sup$-norm with
contractivity constant $|s|$. Hence, $\mathscr{B}$ has a unique fixed point $f\in C(\Delta)$, called a {\em fractal surface over the foldable figure
$\Delta$.} As before, there exists a linear isomorphism $\boldsymbol{\lambda} := (\lambda_1,\ldots,\lambda_N) \mapsto f_{\boldsymbol{\lambda}}$
expressing the fact that the fractal surface $f$ is uniquely determined by the vector of mappings $\lambda_i$.
\begin{example}\label{ex2}
Take the foldable figure on the right-hand side of Figure~\ref{fig5} as the domain $\Delta$ for an affine fractal function. The four subsimplices
$\Delta_1,\ldots, \Delta_4$ induce four similitudes
\begin{align*}
u_1 (x,y) &= \frac{1}{2}\,\begin{pmatrix} 1 & 0\\0 & 1\end{pmatrix}\begin{pmatrix}x \\y\end{pmatrix} + \begin{pmatrix}\frac{1}{2} \\0\end{pmatrix},
& u_2(x,y) &= \frac{1}{2}\,\begin{pmatrix} -1 & 0\\0 & 1\end{pmatrix}\begin{pmatrix}x \\y\end{pmatrix} + \begin{pmatrix}\frac{1}{2}\\0\end{pmatrix}\\
u_3 (x,y) &= \frac{1}{2}\,\begin{pmatrix} 1 & 0\\0 & -1\end{pmatrix}\begin{pmatrix}x \\y\end{pmatrix} + \begin{pmatrix}0 \\ \frac{1}{2} \end{pmatrix},
& u_4(x,y) &= \frac{1}{2}\,\begin{pmatrix} 1 & 0\\0 & 1\end{pmatrix}\begin{pmatrix}x \\y\end{pmatrix} + \begin{pmatrix}0\\ \frac{1}{2}\end{pmatrix}.
\end{align*}
The similarity ratio $\sigma$ equals 1/2. Choose $s := 3/5$ and as functions $\lambda_1,\ldots, \lambda_4$:
\begin{align*}
\lambda_1 (x,y) &:= -\frac{1}{5}\,x + \frac{3}{10}\,y + \frac{1}{5} =: \lambda_2 (x,y),\\
\lambda_3 (x,y) &:= \frac{1}{5}\,x - \frac{3}{10}\,y + \frac{3}{10} =: \lambda_4 (x,y).\\
\end{align*}
A short computation shows that this collection of functions satisfies condition (*). Figure~\ref{fig6} shows the graph of the affine fractal surface
generated by these maps.
\begin{figure}[h]
\begin{center}\includegraphics[width=2.5in,height=1.5in]{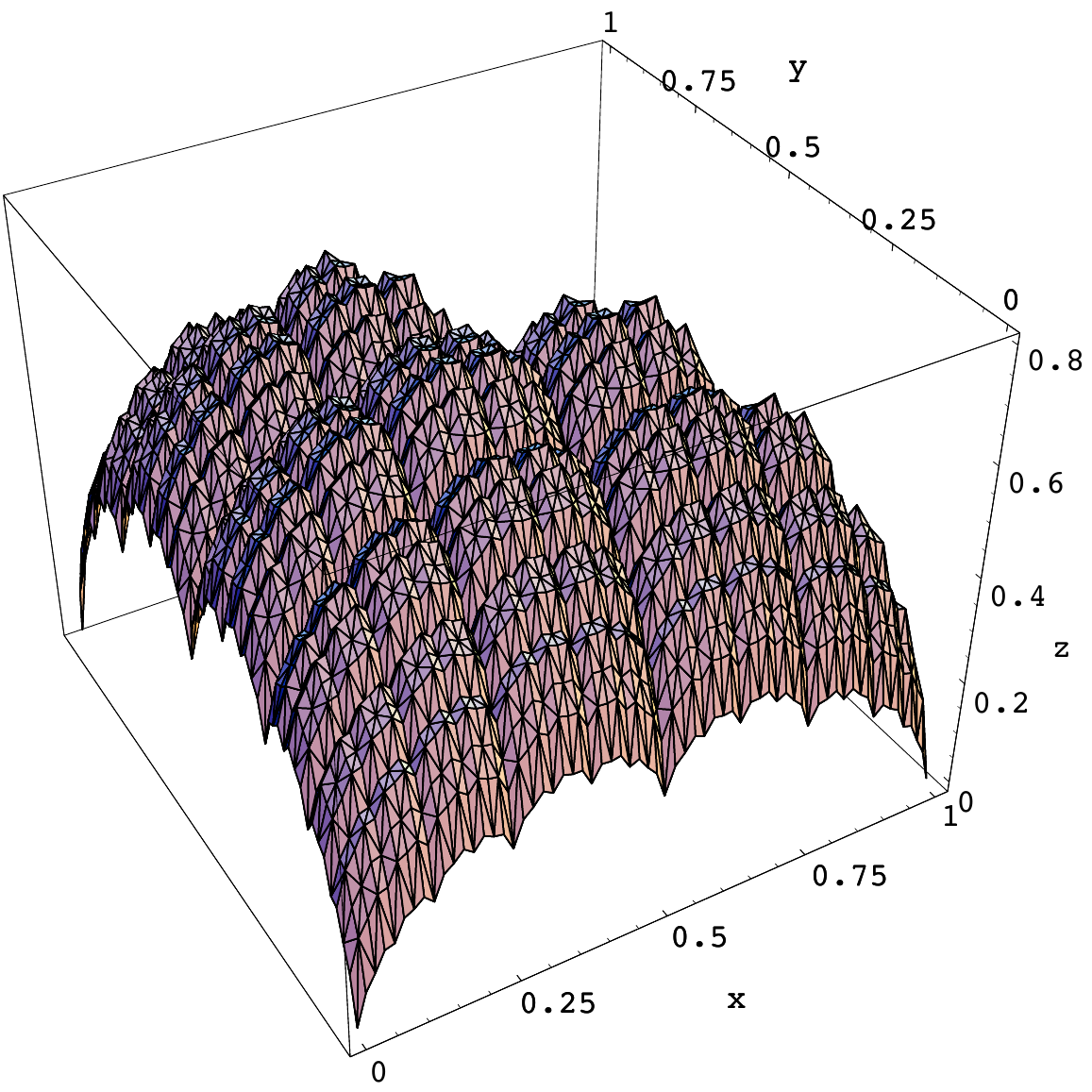}\end{center}
\caption{An affine fractal surface.}\label{fig6}
\end{figure}
\newline
Notice that the affine fractal surface has the value $z = 0$ at the outer vertices $(0,0)$, $(1,0)$, and $(0,1)$ of $\Delta$, and the values
$z = 1/2$, $z = 1/2$, and $z = 3/10$ at the inner vertices $(1/2,0)$, $(1/2,1/2)$, and $(0,1/2)$ of $\Delta$. It is not hard to see that the space
of affine fractal surfaces over $\Delta$ is six-dimensional; there is one basis fractal surface for each outer and inner vertex of $\Delta$. If we
denote the set of inner and outer vertices of $\Delta$ by $\{(x_i,y_j)\st 2 \leq i+j \leq 4\}$ and by $z_{ij}$ the associated $z$ value of $f$, then
\[
f = \sum_{i,j} z_{ij}\,\varphi_{ij},
\]
where $\varphi_{ij}(x_k,y_\ell) = \delta_{ik}\,\delta_{j\ell}$ is a basis fractal surface. Three of these basis fractal surfaces, namely
$\varphi_{12}$, $\varphi_{21}$, and $\varphi_{22}$ are displayed in Figure~\ref{fig7}.
\begin{figure}
\vskip -0.25cm
\includegraphics[width=2in,height=1.5in]{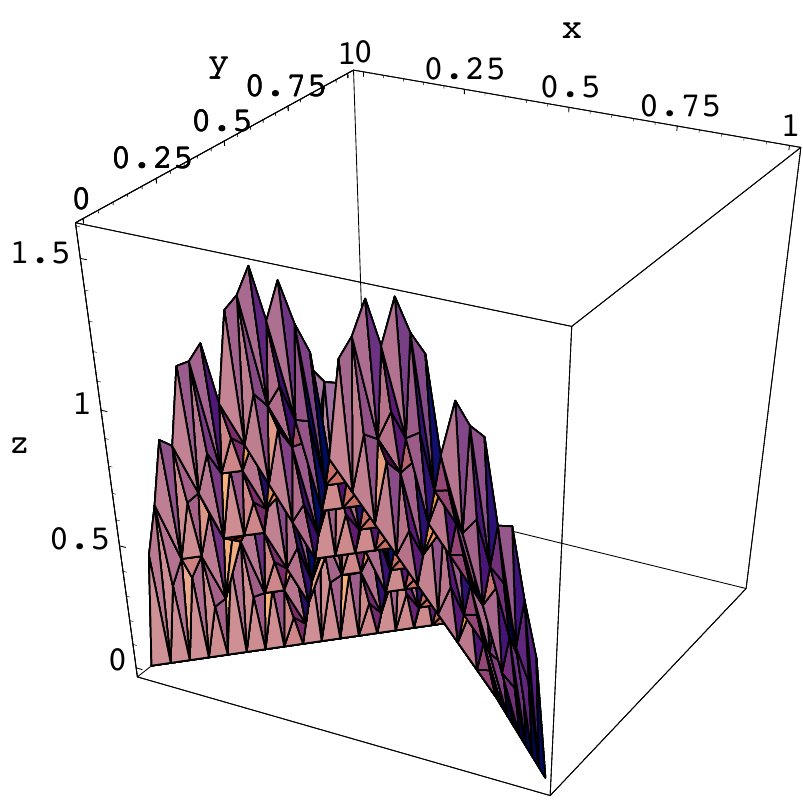}\hspace{1cm}
\includegraphics[width=2in,height=1.5in]{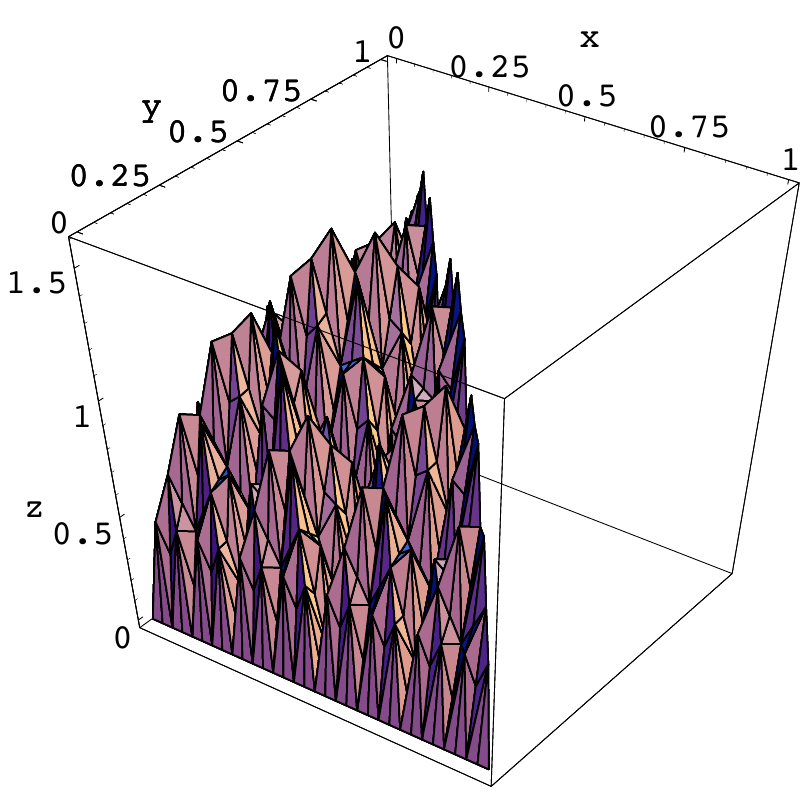}%
\begin{center}\includegraphics[width=2in,height=1.5in]{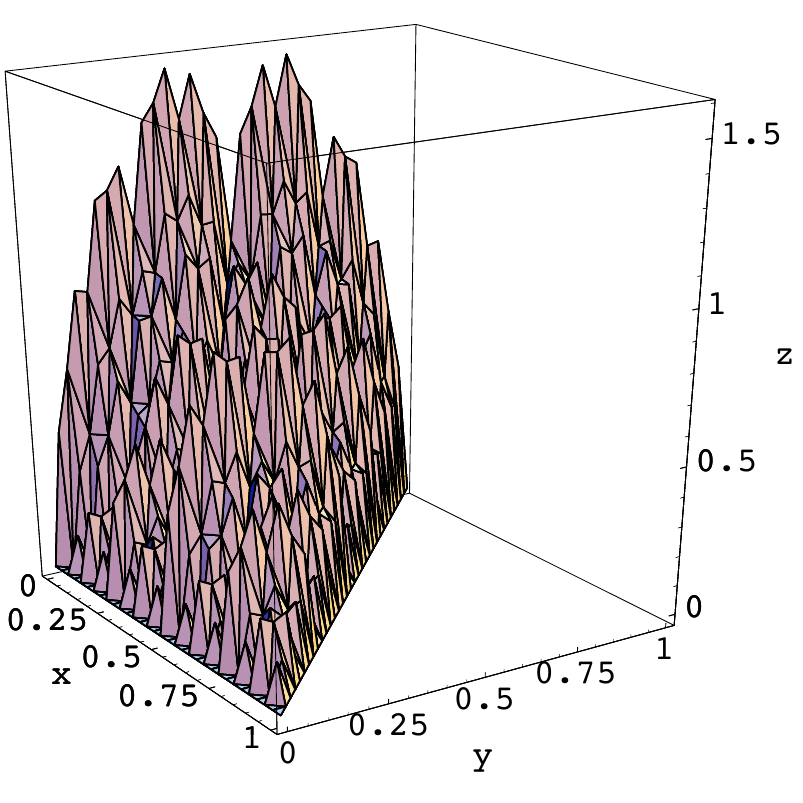}\end{center}
\caption{Three basis fractal surfaces.}\label{fig7}
\end{figure}
\end{example}
\begin{remark}\label{rem1}
Let ${\bf i} := (i_1,\ldots, i_k)\in\{1,\ldots, N\}^k$. Then, if $\{\varphi_\nu\}$ is a fractal function basis for a fractal surface $f$ defined on
a foldable figure $\Delta\in\mathbb{R}^n$, then $\{u_{{\bf i}}^\sharp\varphi_\nu\}$, where $(u_{{\bf i}}^\sharp\varphi_\nu) (x) := \lambda_\nu\circ
u_{\bf i}^{-1} (x) + s_\nu \varphi_\nu\circ u_{\bf i}^{-1}(x)$, $x\in u_{\bf i} (\Delta)$, is a basis for $f\vert_{u_{\bf i} (\Delta)}$.
\end{remark}
\begin{remark}
A more general class of fractal surfaces can be defined by taking the $\lambda_i$'s to be continuous functions satisfying condition (*) and not
merely affine functions. (\cite{GHM1,GHM2,M2})
\end{remark}
\par
In order to achieve out goal, namely to construct a multiresolution analysis on $L^2(\mathbb{R}^n)$ with an orthonormal basis, we need to take into
account the algebraic structure of a foldable figure, i.e., its associated affine Weyl group $\widetilde{\mathcal{W}}$, when constructing affine
fractal surfaces on all of $\mathbb{R}^n$.
\par
To this end, let $F\subset\mathbb{R}^n$ be a foldable figure with $0\in\mathbb{R}^n$ as one of its vertices. Denote by $\mathcal{H}$ be the set of
hyperplanes associated with $F$ and by $\Sigma$ be the tessellation of $F$ induced by $\mathcal{H}$. The affine Weyl group of the foldable figure
$F$, $\widetilde{\mathcal{W}}$, is then the group generated by $\mathcal{H}$. The following theorem summarizes some of the properties of $F$,
$\mathcal{H}$, and $\widetilde{\mathcal{W}}$.
\begin{theorem}
Let $F$ be a foldable figure with associated set of hyperplanes $\mathcal{H}$ and affine Weyl group $\widetilde{\mathcal{W}}$. Then
\begin{enumerate}
\item   $\mathcal{H}$ consists of the translates of a finite set of linear hyperplanes.
\item   $\widetilde{\mathcal{W}}$ is simply-transitive on $\Sigma$, i.e., for all $\sigma,\tau\in\Sigma$ there exists a unique $r\in
        \widetilde{\mathcal{W}}$ such that $\tau = r \sigma$.
\item   $\forall\varkappa\in\mathbb{N}: \varkappa\mathcal{H}\subset\mathcal{H}$. [Here $\varkappa\mathcal{H} := \{\varkappa H\st
        H\in\mathcal{H}\}$]
\end{enumerate}
\end{theorem}
Now, fix $1 < \varkappa\in\mathbb{N}$ and define $\Delta := \varkappa F$. Then $\Delta$ is also a foldable figure, whose $N := \varkappa^n$
subfigures $\Delta_i\in\Sigma$. Assume w.l.o.g that $\Delta_1 = F$. The tessellation and set of hyperplanes for $\Delta$ are $\varkappa \Sigma$ and
$\varkappa\mathcal{H}$, respectively. Moreover, the affine reflection group generated by $\varkappa\mathcal{H}$ is an isomorphic subgroup of
$\widetilde{\mathcal{W}}$. Note that the similarity ratio $\sigma = 1/\varkappa$. By simple transitivity of $\widetilde{\mathcal{W}}$, define
similitudes $u_i:\Delta\to\Delta_i$ by:
\[
u_1 := (1/\varkappa)(\cdot)\quad\text{and}\quad\forall j = 2,\ldots, N:\,\quad u_j := r_{j,1}\circ u_1.
\]
Denote by $\Pi^d = \Pi^d (\mathbb{R}^n)$ the linear space of real polynomials of degree at most $d$ and choose functions $\lambda_1,\ldots, \lambda_N
\in \Pi^d$ satisfying condition $(*)$.\footnote{More general functions can be used, but for the purposes of this paper the restriction to polynomials
provides a large and important subclass.} We denote the linear space of all such polynomials by $J^d$. If we define an operator $\mathscr{B}$ as in
Equation~(\ref{10}) but with $\boldsymbol{\lambda} = (\lambda_1,\ldots,\lambda_N)\in (J^d)^N$, then $\mathscr{B}$ is again contractive in the
$\sup$-norm on $C(\Delta)$ with contractivity $|s|$. Hence, it has a unique fixed point, which is the by $\boldsymbol{\lambda}\in (J^d)^N$ generated
fractal surface $f = f_{\boldsymbol{\lambda}}$ on $\Delta$.
\par
To extend $f$ to all of $\mathbb{R}^n$ we use that fact that the foldable figure $F$ is a fundamental domain for its associated affine Weyl group
and that it tessellates $\mathbb{R}^n$ by reflections in its bounding hyperplanes, i.e., under the action of $\widetilde{\mathcal{W}}$. To this end,
let
\[
J^{\widetilde{\mathcal{W}}}:= \prod \biggl\{(J^d)^N \st r\in\widetilde{\mathcal{W}}\biggr\}.
\]
For $\boldsymbol{\Lambda}\in J^{\widetilde{\mathcal{W}}}$, define $f_{\boldsymbol{\Lambda}}$ by
\[
f_{\boldsymbol{\Lambda}}\vert_{r( \overset{\circ}{\Delta})} := f_{\boldsymbol{\Lambda}(r)}\circ r^{-1},\quad r\in\widetilde{\mathcal{W}},
\]
where $\boldsymbol{\Lambda}(r) = (\boldsymbol{\Lambda}(r)_1,\ldots,\boldsymbol{\Lambda}(r)_N)$ is the $r$-th coordinate of $\boldsymbol{\Lambda}$.
\begin{remark}
The values of $f_{\boldsymbol{\Lambda}}$ are left unspecified on the hyperplanes $\varkappa\mathcal{H}$, a set of Lebesgue measure zero in
$\mathbb{R}^n$, and thus $f_{\boldsymbol{\Lambda}}$ actually represents an equivalence class of functions.
\end{remark}
\section{Dilation- and $\widetilde{\mathcal{W}}$-invariant spaces}
Let $V$ be a linear space of functions  $f: \mathbb{R}^n\to\mathbb{R}$, $D$ an expansive unitary operator on $V$, and $\widetilde{\mathcal{W}}$ an affine Weyl group. (By an expansive unitary operator on $V$ we mean a dilation operator of the form given in (\ref{e1}) for an expansive matrix on $\R^n$.)
\begin{definition}
The linear space $V$ is called {\em dilation--invariant} if
\[
D^{-1} V \subset V
\]
and $\widetilde{\mathcal{W}}$--invariant if
\[
f\in V \Longrightarrow f\circ r \in V,\quad\forall r\in\widetilde{\mathcal{W}}.
\]
\end{definition}
Dilation--invariance of a global fractal function $f_{\boldsymbol{\Lambda}}$ can be expressed in terms of an associated dilation--invariance of $\boldsymbol{\Lambda}\in J^{\widetilde{\mathcal{W}}}$.
\begin{theorem}
Let $1 < \varkappa \in\mathbb{N}$, $D_\varkappa := \varkappa\,\textrm{id}_{\mathbb{R}^n}$ and $f_{\boldsymbol{\Lambda}}$ a global fractal function
generated by $\boldsymbol{\Lambda}\in J^{\widetilde{\mathcal{W}}}$. Then
\begin{enumerate}
\item   $f_{\boldsymbol{\Lambda}}$ is $D_\varkappa$--invariant on $\mathbb{R}^n$ iff $\boldsymbol{\Lambda}$ is $\delta_\varkappa$--invariant
        on $J^{\widetilde{\mathcal{W}}}$, in other words, iff
        \[
        D_\varkappa^{-1} f_{\boldsymbol{\Lambda}} = f_{\delta_\varkappa \boldsymbol{\Lambda}},
        \]
        where $\delta_\varkappa: J^{\widetilde{\mathcal{W}}}\to J^{\widetilde{\mathcal{W}}}$, is given by
        \[
        \hspace*{1.25cm}\delta_\varkappa \boldsymbol{\Lambda} (\varkappa\,r\,u_j)_i = \boldsymbol{\Lambda}(r)_j\circ u_i + s\,[\boldsymbol{\Lambda}
        (r)_i - \boldsymbol{\Lambda}(r)_j],\quad r\in\widetilde{\mathcal{W}},\, i,j = 1, \ldots, N.
        \]
\item   $J^{\widetilde{\mathcal{W}}}$ is $\delta_\varkappa$--invariant.
\end{enumerate}
\end{theorem}
\begin{proof}
For part (1) see \cite{GHM2,M2}, part (2) is a simple calculation.
\end{proof}
The above theorem now allows us to generate a multiresolution analysis on $L^2(\mathbb{R}^n)$ and define an orthonormal basis for it consisting of
fractal surfaces generated by functions in $J^{\widetilde{\mathcal{W}}}$.
\begin{definition}
A multiresolution analysis (MRA) of $L^2(\mathbb{R}^n)$ with respect to dilation $D_\varkappa$ and affine Weyl group $\widetilde{\mathcal{W}}$ consists of a sequence of spaces $\{V_k\st k\in\mathbb{Z}\}\subset L^2(\mathbb{R}^n)$ satisfying
\begin{enumerate}
\item   $V_k\subset V_{k+1}$, for all $k\in\mathbb{Z}$.
\item   $\bigcup \{V_k\st k\in\mathbb{Z}\}$ is dense in $L^2(\mathbb{R}^n)$.
\item   There exists a finite set of generators $\{\phi^a\st a\in A\}\subset L^2(\mathbb{R}^n)$ such that
        \[
        \mathcal{B}_\phi := \{\phi^a\circ r\st a\in A,\,r\in\widetilde{\mathcal{W}}\}
        \]
        is a Riesz basis for $V_0$.
\item   $f\in V_k$ $\Longrightarrow$ $D_\varkappa f\in V_{k+1}$, $\forall k\in\mathbb{Z}$.
\end{enumerate}
\end{definition}
\noindent
We know that $J^{\widetilde{\mathcal{W}}}$ is $\delta_\varkappa$--invariant and thus we define
\[
V_0 := \left\{f_{\boldsymbol{\Lambda}}\st\boldsymbol{\Lambda}\in J^{\widetilde{\mathcal{W}}}\right\}
\]
and
\[
V_k := D_\varkappa^{k} V_0, \qquad\forall k\in\mathbb{N}.
\]
Notice that $D_\varkappa$ can be represented by the expansive matrix $\varkappa I$, where $I$ is the $n\times n$-identity matrix. In the following,
we will not distinguish between the operator $D_\varkappa$ and its matrix representation and denote both by $D_\varkappa$.
\par
Since the dimension of $(J^d)^N = (d+1)N$, the dimension of $\dim V_0\vert_\Delta$ is also $(d+1)N$, and thus $|A| = (d+1) N$. We can take as a
basis,a fractal surface basis of the type considered in Example~\ref{ex2}, and apply the Gram-Schmidt orthogonalization procedure to it, to obtain an
orthonormal basis $\{\phi^a\st a = 1,\ldots, |A|\}$ for $V_0\vert_\Delta$.
\par
Let $\Phi := (\phi^1,\ldots, \varphi^{|A|})^\top$. Then $V_1 \subset V_0$ implies the existence of sequence of $|A| \times |A|$-matrices $\{P(r)\st
r\in\widetilde{\mathcal{W}}\}$, only a finite number of which are nonzero, such that
\be\label{ref}
(D_\varkappa ^{-1}\Phi)(x) = \Phi (x/\varkappa) = \sum_{r\in\widetilde{\mathcal{W}}} P(r)\,(\Phi \circ r)(x).
\ee
Equation~(\ref{ref}) is the refinement equation for the scaling vector $\Phi$ in the current situation.
\begin{theorem}
The ladder of spaces $\{V_k\st k\in\mathbb{Z}\}$ defines an MRA of $L^2(\mathbb{R}^n)$ with respect to dilation $D_\varkappa$ and affine Weyl group
$\widetilde{\mathcal{W}}$.
\end{theorem}
\begin{proof}
See \cite{GHM1,GHM2,M2}.
\end{proof}
\par
For $k\in\mathbb{N}$, define the wavelet spaces $W_k := V_{k+1} \ominus V_k$. Since $\dim W_0\vert_\Delta = \dim V_1\vert_\Delta - \dim V_0\vert_
\Delta = (\varkappa^n - 1)(d+1)N$, we can again use the Gram-Schmidt orthonormalization procedure to construct an orthonormal basis $\{\psi^b\st
b \in B\}$ for $W_0\vert_\Delta$, where $B$ has cardinality $|B|= (\varkappa^n - 1)(d+1)N$.
\par
Let $\Psi := (\psi^1,\ldots, \psi^{|B|})^\top$. As $W_0 \subset V_1$ there exists a sequence of $|B| \times |A|$-matrices $\{Q(r)\st r\in
\widetilde{\mathcal{W}}\}$, only a finite number of which are nonzero, such that
\be
(D_\varkappa ^{-1}\Psi)(x) = \Psi (x/\varkappa) = \sum_{r\in\widetilde{\mathcal{W}}} Q(r)\,(\Phi \circ r)(x).
\ee
In addition, one can write down finite decomposition and reconstruction algorithms for these multigenerators $\Phi$ and $\Psi$. The interested reader
is referred to \cite{GHM1,GHM2,M2} for further details.
\section{Wavelet sets constructed via Coxeter groups}
In this section, a new type of wavelet set is introduced which we
will call a \emph{dilation-reflection wavelet set}. It belongs to
the Coxeter/fractal-surface multiresolution analysis theory.  The
idea is to adapt Definition \ref{dil-trans}, replacing the group of translations
$\mathcal{T}$ in the traditional wavelet theory by an affine Weyl
group whose fundamental domain is a foldable figure $C$, and to use
the orthonormal basis of fractal surfaces constructed in the
previous section.
\begin{remark}
In sections 3 - 6, the requirement that 0 is a vertex of $C$ is not necessary. We chose it for the construction of fractal surfaces since the maps
are easier to define (otherwise a shift is to be added) and this was consistent with earlier treatments of the subject \cite{GHM1,GHM2,M2}. In the
present section, however, where we define a dilation-reflection wavelet set, there is a disadvantage in requiring that 0 is a vertex of $C$. It
agrees more with the dilation theory if 0 is an interior point of $C$ to simplify the application of Definition \ref{dil-trans} and Theorem \ref{t1}
to this setting to produce the dilation-reflection wavelet sets. This requires an affine shift (as mentioned above) in the Coxeter/fractal surface
multiresolution theory discussed in the previous sections. The expositional stance we take is to give a formal definition (Definition \ref{def7.4})
of dilation-reflection wavelet set in more general terms involving an affine shift in both the Weyl (abstract translation) group and the matricial
(abstract dilation) group. We simply take the dilation group fixed point in Definition \ref{dil-trans} to be any point $\theta$ in the (nonempty)
interior of $C$. We leave the details involving the affine shift in the Weyl group to the reader because they are straight-forward, and give explicit
details on how the dilation group needs to be affinely shifted to agree with Definition \ref{dil-trans} and Theorem \ref{t1}. For the concrete
examples \ref{ex1} and \ref{ex2}, the theory is the simplest from the dilation group viewpoint in the case where 0 is in the interior of $C$ and we
therefore take $\theta = 0$. Hopefully, these concrete examples will clarify our treatment of the theory.
\end{remark}
In Definition \ref{dil-trans}, we take $X := \mathbb{R}^n$ endowed with the Euclidean affine structure and distance, and for the abstract translation group $\mathcal{T}$ we take the affine Weyl group $\widetilde{{\mathcal{W}}}$ generated by a group of affine reflections arising from a locally finite collection of affine hyperplanes of $X$. Let $C$ denote a fundamental domain for $\widetilde{{\mathcal{W}}}$ which is also a foldable figure. Recall that $C$ is a simplex\footnote{Let $\{x_0,x_1,\ldots, x_n\}$ be a set of linearly independent points in $\R^n$. The set 
\[
\Sigma^n := \left\{\sum_{i=0}^n \lambda_i\,x_i \Bigg\vert \sum_{i=0}^n \lambda_i =1, \; \lambda_i \geq 0,\, i = 0,1,\ldots n\right\}
\]
is called a simplex.}, i.e., a convex connected polytope (here we do not assume it has $n+1$ vectors), which tessellates $\mathbb{R}^n$ by reflections about its bounding hyperplanes. Let $\theta$ be any fixed interior point of $C$. Let $A$ be any real expansive matrix in $M_n(\mathbb{R})$ acting as a linear transformation on $\mathbb{R}^n$.  In the case where $\theta$ is the orgin $0$ in $\mathbb{R}^n$ we simply take $D$ to be the usual dilation by $A$ and the abstract dilation group to be $\mathcal{D} = \{D^k \st k \in \mathbb{Z}\}$. For a general $\theta$, define $D$ to be the affine mapping $D(x) := A(x - \theta) + \theta, x \in\mathbb{R}^n$ and $\mathcal{D}_{\theta} = \{D^k \st k \in \mathbb{Z}\}$.
\begin{proposition}\label{th10}
 $(\mathcal{D}_\theta,\widetilde{{\mathcal{W}}})$ is an
abstract dilation-translation pair in the sense of Definition \ref{dil-trans}.
\end{proposition}
\begin{proof}
By the definition of $D$, $\theta$ is a fixed point for $\mathcal{D}_\theta$. By a change of coordinates we may assume without loss of generality that $\theta = 0$ and consequently that $D$ is multiplication by $A$ on $\mathbb{R}^n$.
\par
Let $B_r (0)$ be an open ball centered at $0$ with radius $r > 0$ containing both $E$ and $C$. Since $F$ is open and $A$ is expansive, there exists a $k \in\N$ sufficiently large so that $D^k F$ contains an open ball $B_{3r}(p)$ of radius $3r$ and with some center $p$.  Since $C$ tiles $\R^n$ under the action of $\wW$, there exists a word $w\in \wW$ such that $w(C) \cap B_r(p)$ has positive measure. (Note here that $B_r(p)$ is the ball with the same center $p$ but with smaller radius $r$.) Then $w(B_r(0)) \cap B_r(p) \neq \emptyset$. Since reflections (and hence words in $\wW$) preserve diameters of sets in $\R^n$, it follows that $w(B_r(0))$ is contained in $B_{3r}(p)$.  Hence $w(E)$ is contained in $D^k(F)$, as required.
\par
This establishes part (1) of Definition \ref{dil-trans}. Part (2) follows from the fact that $\theta = 0$ and $D$ is multiplication by an expansive matrix in $M_n(\mathbb{R})$.
\end{proof}
\begin{definition} \label{def7.2}
Given an affine Weyl group $\widetilde{{\mathcal{W}}}$ acting on $\mathbb{R}^n$ with fundamental domain a foldable figure $C$, given a designated interior point $\theta$ of $C$, and given an expansive matrix $A$ on $\mathbb{R}^n$, a dilation--reflection wavelet set for $(\widetilde{{\mathcal{W}}}, \theta, A)$ is a measurable subset $E$ of $\mathbb{R}^n$ satisfying the properties:
\begin{enumerate}
\item $E$ is congruent to $C$ (in the sense of Definition 2.4) under the
action of $\widetilde{{\mathcal{W}}}$, and
\item $W$ generates a measurable partition of $\mathbb{R}^n$ under the action of
the affine mapping $D(x) := A(x - \theta) + \theta$.
\end{enumerate}
In the case  where $\theta = 0$, we abbreviate $(\widetilde{{\mathcal{W}}}, \theta, A)$ to $(\widetilde{{\mathcal{W}}}, A)$.
\end{definition}
\begin{theorem}\label{100}
There exist $(\widetilde{{\mathcal{W}}}, \theta, A)$--wavelet sets for every choice of $\widetilde{W}$, $\theta$, and $A$.
\end{theorem}
\begin{proof}
This is a direct application of Theorem \ref{t1}. Let $C$ be a fundamental domain for $\widetilde{W}$ which is a foldable figure, let $\theta$ be any interior point of $C$, and let $A$ be any expansive matrix in $M_n (\mathbb{R})$. By Proposition \ref{th10}, $(\mathcal{D}_\theta, \widetilde{{\mathcal{W}}})$ is an abstract dilation--translation pair with $\theta$ the dilation fixed point. Let $C$ play the role of $E$ in Theorem \ref{t1}. As in the proof of Theorem \ref{Th2.11}, which is sketched above the statement of the theorem, let $F_A := A(B)\setminus B$, where $B$ is the unit ball of $\mathbb{R}^n$, and let $F:= F_A + \theta$. Then  $\{D^k F\st k \in \mathbb{Z}\}$ is a partition of $\mathbb{R}^n\setminus\theta$, where $D$ is the affine map $D(x) := A(x_\theta) + \theta$. Since $F$ has nonempty interior and is bounded away from $\theta$, Theorem \ref{t1} applies yielding a measurable set $W$ which is simultaneously congruent to $C$ under the action of $\widetilde{W}$ and congruent ot $F$ under the action of $\mathcal{D}_\theta$. Since $F$ generates a measurable partition of $\mathbb{R}^n$ under $\mathcal{D}_\theta$, so must any set that is $\mathcal{D}_\theta$-congruent to $F$. Hence $W$ satisfies (2) of Definition \ref{def7.2}. Since it is also $\widetilde{W}$--congruent to $C$, this shows that it is a dilation--reflection wavelet set for $(\widetilde{{\mathcal{W}}}, \theta, A)$, as required.
\end{proof}
\begin{remark}\label{rem7.5}
The role of $C$ in the dilation--reflection wavelet theory is
analogous to the role of the interval $[0, 2\pi)$ in the dyadic
dilation--translation wavelet theory on the real line. For sake of
exposition, let us recapture this role:  The set of exponentials
\begin{equation}
\left\{\frac{e^{i\ell s}}{\sqrt{2\pi}}\Big|_{[0, 2\pi)}\Bigg\vert\;\ell\in
\mathbb{Z}\right\}
\end{equation}
is an orthonormal basis for $L^2([0, 2\pi))$, hence if $W$ is any set
which is $2\pi$--translation congruent to $[0,2\pi)$, then
\begin{equation}
\left\{\frac{e^{i\ell s}}{\sqrt{2\pi}}\Big|_{W}\Bigg\vert\;\ell\in
\mathbb{Z}\right\}
\end{equation}
is an orthonormal basis for $L^2(W)$.  A dyadic dilation--translation
wavelet set on the line has this (spectral set) property, and also
generates a measurable partition of $\mathbb{R}$ under dilation by 2, and
consequently the union of the sets
\begin{equation}
D^n\left\{\frac{e^{i\ell s}}{\sqrt{2\pi}}\Big|_{W}\Bigg\vert\; \ell\in
\mathbb{Z}\right\}
\end{equation}
is an orthonormal basis for $L^2(\mathbb{R})$.
\par
So, recapitulating, the role of $[0, 2\pi)$ is that it supports a
``special" orthonormal basis for $L^2([0, 2\pi))$ induced by the
translation group via the Fourier Transform, and thus $W$, being
$\tau$--congruent to $[0, 2\pi)$, also supports an orthonormal basis
for $L^2(W)$ induced by the $\tau$--congruence. The role of the
fundamental domain $C$ in the dilation-reflection theory is analogous
to this.
\end{remark}
For sake of exposition, it is natural to make the following somewhat
abstract definition.
\begin{definition}\label{def7.4}
An abstract wavelet set in $\mathbb{R}^n$ is a measurable set
$W$ that produces an orthonormal basis for $L^2(\mathbb{R}^n)$ under the action
of two countable unitary systems acting consecutively, with the
first system inducing an orthonormal basis for $L^2(W)$ by its
action restricted to $W$, and the second is a system of dilations by a
family of affine transformations on $\mathbb{R}^n$ whose action on $W$ yield a
measurable partition of $\mathbb{R}^n$ (and hence the dilation unitaries when
applied to $L^2(W)$ yield a direct-sum orthogonal decomposition of
$L^2(\mathbb{R}^n)$.
\end{definition}
In the case of a dilation--translation wavelet set $W$, the two systems of unitaries are $\mathcal{D} := \{D_A^k \st k\in\mathbb{Z}\}$, where $A\in M_n (\mathbb{R})$ is an expansive matrix, and $\mathcal{T} := \{T^\ell \st \ell\in\mathbb{Z}^n\}$. An orthonormal wavelet basis of $L^2 (\mathbb{R}^n)$ is then obtained by setting $\widehat{\psi}_W := (m(W))^{-1/2} \chi_W$ and taking
\[
\left\{\widehat{D}_A^k \widehat{T}^\ell \widehat{\psi}_W \st k\in\mathbb{Z}, \ell\in\mathbb{Z}^n\right\}.
\]
\par
For the systems of unitaries  $\mathcal{D} := \{D_A^k \st k\in\mathbb{Z}\}$ and $\widetilde{\mathcal{W}}$, the affine Weyl group associated with a foldable figure $C$, one obtains as an orthonormal basis for $L^2 (\mathbb{R}^n)$
\[
\left\{D_\varkappa^k \mathcal{B}_\phi \st k\in\mathbb{Z}\right\},
\]
where  $\mathcal{B}_\phi = \left\{\phi^a\circ r\st a\in A,\,r\in\widetilde{\mathcal{W}}\right\}$ is a fractal surface basis as constructed in the
previous section.
\vskip 5pt
In \cite{DL}, two examples of wavelet sets are given in the plane for dilation--2 (i.e., the dilation matrix is $A := 2 I$, where $I$ is the
$2 \times 2$ identity matrix) and $2\pi$-translation (separately in each coordinate). Both examples are reproduced here and it is shown that the
two dilation--translation wavelet sets are also dilation-reflection wavelet sets for $(\widetilde{W}, \theta, 2 I)$, where $\widetilde{W}$ is the
affine Weyl group generated by the reflections about the bounding hyperplanes of $C := [-\pi, \pi)\times [-\pi,\pi)$, and where $\theta = 0$. In
other words, we take dilation to be exactly the same and replace the $2\pi$--translation by the action of the Weyl group. We find it quite
interesting that the same measurable set is a wavelet set in each of the two different theories because the Weyl group and the usual translation
by $2\pi$--group act completely differently as groups of transformations of $\mathbb{R}^2$. In some sense, this justifies our usage of the term
``wavelet set" to denote our construction in the Coxeter/fractal surface theory. The other reason is the interpretation given in Remark \ref{rem7.5}.
\begin{example}\label{ex1}
Let $A := 2 I$, where $I$ denotes the  identity matrix in $\mathbb{R}^2$. For $n\in\mathbb{N}$, define vectors $\vec{\alpha},\vec{\beta}\in\mathbb{R}^2$ by
\[
\begin{split}
\vec{\alpha}_n &:= \frac{1}{2^{2n-2}}\begin{pmatrix}\frac{\pi}{2}\\ \frac{\pi}{2}\end{pmatrix};\\
\vec{\beta}_0 & := 0, \quad \vec{\beta}_n := \sum_{k=1}^n \vec{\alpha_k}.
\end{split}
\]
Define
\[
\begin{split}
G_0 &:= \left[0,\frac{\pi}{2}\right];\\
G_n &:= \frac{1}{2^{2n}}G_0 + \vec{\beta}_n;\\
E_1 &:= \bigcup_{k=1}^\infty G_k \subset 2G_0\setminus G_0;\\
C_1 &:= G_0\cup E_1 + \begin{pmatrix} 2\pi\\ 2\pi\end{pmatrix};\\
B_1 &:= 2G_0\setminus (G_0\cup E_1).
\end{split}
\]
Finally, let
\[
\begin{split}
A_1 &:= B_1\cup C_1;\\
A_2 &:= \{(-x,y) \st (x,y)\in A_1\};\\
A_3 &:= \{(-x,-y) \st (x,y)\in A_1\};\\
A_4 &:= \{(x,-y) \st (x,y)\in A_1\};\\
W_1 &:= A_1\cup A_2\cup A_3\cup A_4.
\end{split}
\]
It is not hard to verify that $W_1$ is $2\pi$-translation--congruent to $C = [-\pi,\pi)\times [-\pi,\pi)$ and a $2$-dilation generator of a measurable partition for the two-dimensional plane $\mathbb{R}^2\setminus \{0\}$. (Cf. \cite{DL})
\par
The set $C$ is a fundamental domain of the (reducible) affine Weyl group $\widetilde{\mathcal{W}}$ that is generated by the (affine) reflections about the six lines, namely $L^x_{\pm}: x = \pm \pi$, $L^y_{\pm}:y = \pm\pi$, $L^x_0: x = 0$, and $L^y_0: y = 0$. (Cf. Left-hand side of Figure \ref{fig5}.) The roots for this Weyl group are given by $r_1  = -r_2 = (1,0)^\top$ and $r_3 = -r_4 = (0,1)^\top$. As a matter of fact, the \textit{Coxeter group} associated with $C$ is Klein's Four Group or the dihedral group $D_4$. Denote by $B_j$, $C_j$, and $E_j$, $j = 2,3,4$, the extension of $B_1$, $C_1$, and $E_1$, respectively, into the $j$th quadrant.  Let $\rho^x_{-}$ and $\rho^y_{-}$ denote the affine reflection about the line $L^x_{-}$ and $L^y_{-}$, respectively. Then it is easily verified that $\rho^y_{-}\rho^x_{-} (E_3)\cup B_1 = 2G_0$. Analogous arguments applied to $E_1$, $E_2$, and $E_4$ show that $W_1$ is $\widetilde{\mathcal{W}}$--congruent to $C$.
\par
Since $C$ is a foldable figure, there exists an orthonormal basis for $L^2 (C)$ generated by fractal surfaces. Now, the sets $2^{-2n}G_0$ are copies of $C$ under appropriate combinations of the maps $u_i$ and we have, by Remark \ref{rem1}, for all $n\in\mathbb{N}$ an orthonormal basis. Hence, the sets $G_n$ have such a basis for all $n\in\mathbb{N}$ and since they only intersect on a set of measure zero, so do the sets $E_j$ and $G_0 \cup E_j$, $j = 1,\ldots,4$. Since the sets $C_j$ are obtained by applying elements of the affine Weyl groups, i.e., isometries, to $G_0 \cup E_j$, $j = 1,\ldots,4$, one obtains an $L^2$-orthonormal basis $\mathcal{B}_\phi$ for $W_1$. Then, since $L^2(W_1)$ is wandering for $D_{2}^k$, $\left\{D_2^k \mathcal{B}_\phi : k\in\mathbb{Z}\right\}$ is an orthonormal basis for $L^2(\mathbb{R}^2)$.
\par
The wavelet set $W_1$ is depicted in Figure \ref{figw1} (See also \cite{DL}).
\begin{figure}[h]
\begin{center}
\includegraphics[width=2in,height=2in]{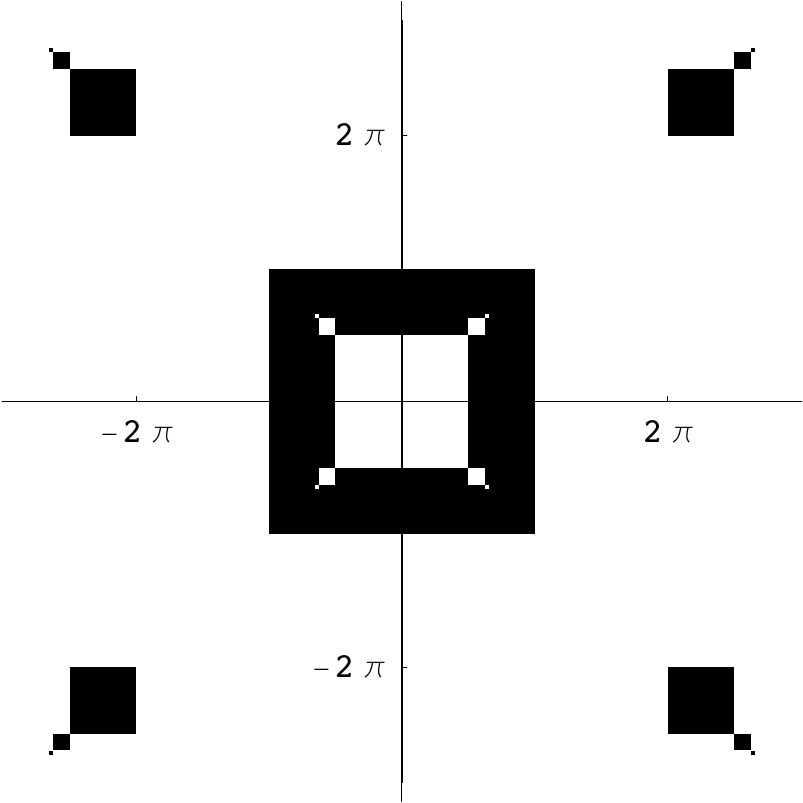}
\caption{The wavelet set $W_1$ in $\mathbb{R}^2$.}\label{figw1}
\end{center}
\end{figure}
\end{example}
\begin{example}\label{example2}
Again, let $A := 2 I$, where $I$ denotes the  identity matrix in $\mathbb{R}^2$. For $n\in\mathbb{N}$, define vectors $\vec{\alpha},\vec{\beta}\in
\mathbb{R}^2$ by
\[
\begin{split}
\vec{\alpha}_n &:= \frac{1}{2^{2n-2}}\begin{pmatrix}\frac{\pi}{2}\\ 0\end{pmatrix};\\
\vec{\beta}_0 & := 0, \quad \vec{\beta}_n := \sum_{k=1}^n \vec{\alpha_k}.
\end{split}
\]
Define
\[
\begin{split}
G_0 &:= \left[0,\frac{\pi}{2}\right]\times \left[-\frac{\pi}{2},\frac{\pi}{2}\right];\\
G_n &:= \frac{1}{2^{2n}}G_0 + \vec{\beta}_n;\\
E &:= \bigcup_{k=1}^\infty G_k \subset 2G_0\setminus G_0;\\
D &:= G_0\cup E + \begin{pmatrix} 2\pi\\ 0\end{pmatrix};\\
B &:= 2G_0\setminus (G_0\cup E).
\end{split}
\]
Define
\[
\begin{split}
A_1 &:= B\cup D;\\
A_2 &:= \{(-x,y) \st (x,y)\in A_1\};\\
W_2 &:= A_1\cup A_2.
\end{split}
\]
That $W_2$ is a dilation--translation wavelet set was established in \cite{DL}.  To show that it is also a dilation--reflection wavelet set, note that $C$ is the same foldable figure as in Example \ref{ex1} above with the same (reducible) affine Weyl group $\widetilde{\mathcal{W}}$. Denote by $L_1: x = {-\pi}$ and $L_2: x = {\pi}$ the left and right bounding lines of $C$ and by $\rho_1$, respectively, $\rho_2$ the corresponding reflections. Let $D^- := \{(-x,y) \st (x,y)\in D\}$. Then it is easy to see that $\rho_2 (D)\cup \rho_1 (D^-) = C$, hence that $W_2$ is $\widetilde{\mathcal{W}}$--congruent to $C$. Arguments similar to those used in Example \ref{ex1} show that there exists an orthonormal basis for $L^2 (\mathbb{R}^2)$ consisting of fractal surfaces.
\par
The wavelet set ${W_2}$ is shown in Figure \ref{figw2}.
\begin{figure}[h]
\begin{center}
\includegraphics[width=7cm,height=2.8cm]{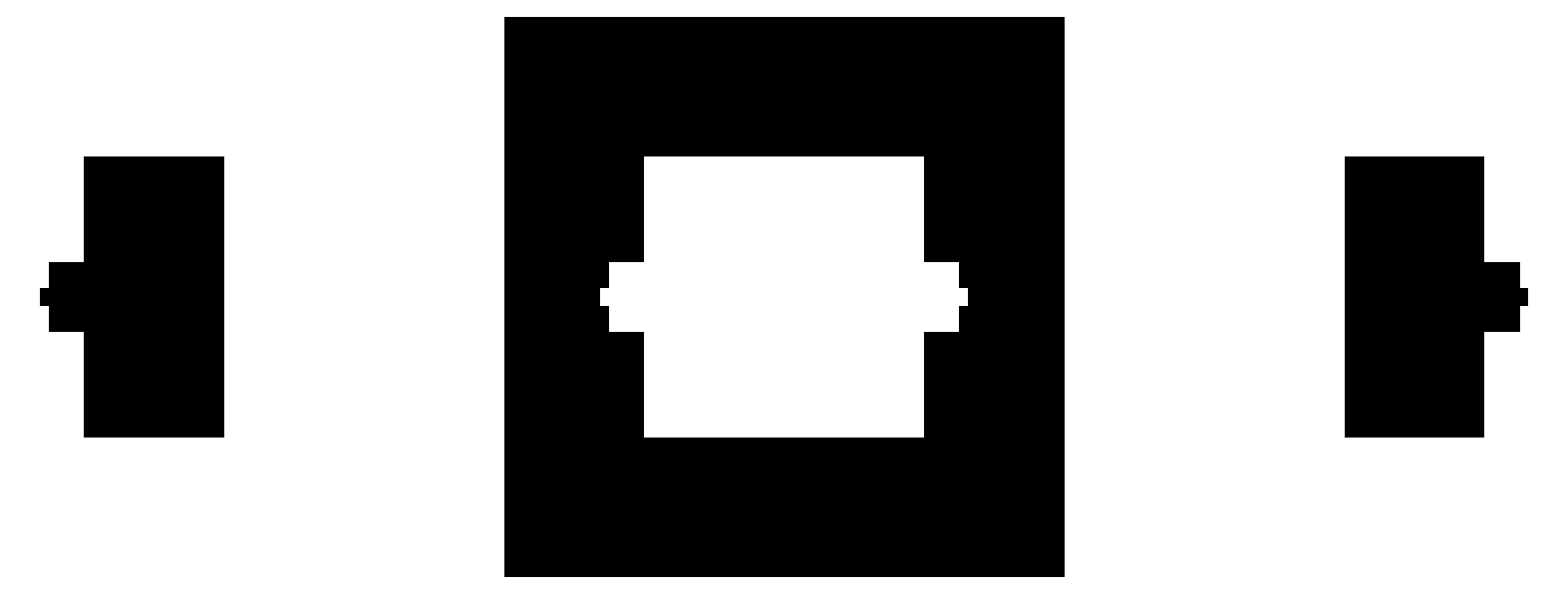}
\caption{The wavelet set $W_2$ in $\mathbb{R}^2$.}\label{figw2}
\end{center}
\end{figure}
\end{example}
\section{Generalitites and Open Problems}
The two examples that were considered at the end of the previous section are representatives of what one might call ``three-way tiling sets" of the Euclidean plane $\mathbb{R}^2$. To this end, observe that if that $\mathcal{G}$ is a group of transformation in $\mathbb{R}^n$ and some subset $K$ tiles $\mathbb{R}^n$ under the action of $\mathcal{G}$, then any set $L$ that is $\mathcal{G}$--congruent to $K$ also tiles $\mathbb{R}^n$ under the action of $\mathcal{G}$. The sets $W_i$, $i=1,2$, are wavelet sets for both the affine reflection group $\widetilde{\mathcal{W}}$ and the translation group $\mathcal{T}$. Both clearly tile $\mathbb{R}^2$ under translation, but they also tile under reflections since both are $\widetilde{\mathcal{W}}$--congruent to a foldable figure, namely $C = [-\pi,\pi]\times [-\pi,\pi]$, which tiles the plane under reflections in its bounding hyperplanes.
\par
The fact that $W_1$ and $W_2$ tile under both the dilation and translation group makes them dilation-translation wavelet sets, and the fact that they tile under both the dilation and reflection group makes them dilation-reflection wavelet sets. Indeed, more can be said about these two wavelets sets and the fundamental domain $C$.
\par
The foldable figure $C$ is in both situations a fundamental domain for the affine Weyl group $\widetilde{\mathcal{W}}$ and for the standard translation-by-$2\pi$ group $\mathcal{T}$. Recall that by Theorem \ref{th4.5}, $\widetilde{\mathcal{W}} = \mathcal{W}\ltimes \Gamma$ where $\mathcal{W}$ is the stabilizer (isotropy group) of the origin and $\Gamma$ the translation group generated by the coroots $r_j^\vee = 2\pi r_j$, $j = 1,\ldots, 4$ (cf. Example \ref{ex1}). (Note that the $k$ in Equation \ref{hyper} is $\pi k$ in the two examples.) If we denote the intersection group $\widetilde{\mathcal{W}}\,\cap\,\mathcal{T}$ by $\mathcal{J}$, then we have in this case $\mathcal{J} \subset \mathcal{T}$. Indeed, $\mathcal{J}$ is generated by translations of $(k\cdot 4\pi, \ell\cdot 4\pi)$, $k, \ell\in\mathbb{Z}$. To see this, we compute $\rho_{r,\pi k}\circ\rho_{s,\pi\ell}$, where $r$ and $s$ are any two of the roots $r_j$, $j=1,\ldots, 4$, and find that for $s = \pm r$
\[
\rho_{r,\pi k}\circ\rho_{\pm r,\pi\ell}  = \textrm{id}_{\mathbb{R}^2} \pm 2\pi (k - \ell) r,\qquad | k - \ell |\geq 2,
\]
whereas for $s\perp r$,
\[
\rho_{r,\pi k}\circ\rho_{s,\pi\ell} = -\textrm{id}_{\mathbb{R}^2} + 2\pi (k r + \ell s).
\]
Thus, every element of $\widetilde{\mathcal{W}}$ is the product of a simple reflection and a translation by  $(k\cdot 4\pi, \ell\cdot 4\pi)$, $k, \ell\in\mathbb{Z}$. Furthermore, $\mathcal{J}$ is clearly nontrivial but also big enough so that it together with the dilation-by-$2 I$ group satisfies the axioms in the definition of ``abstract dilation-translation pair" (Definition \ref{dil-trans}). In addition, both sets $W_i$, $i =1,2$, are actually congruent to $C$ via the intersection group $\mathcal{J}$. In other words, for the translation-by-$2\pi$ congruence only translations in $\mathcal{J}$ are used, and likewise for the $\widetilde{\mathcal{W}}$--congruence.
\par
\begin{figure}[h]
\begin{center}
\includegraphics[width=3cm,height=2.5cm]{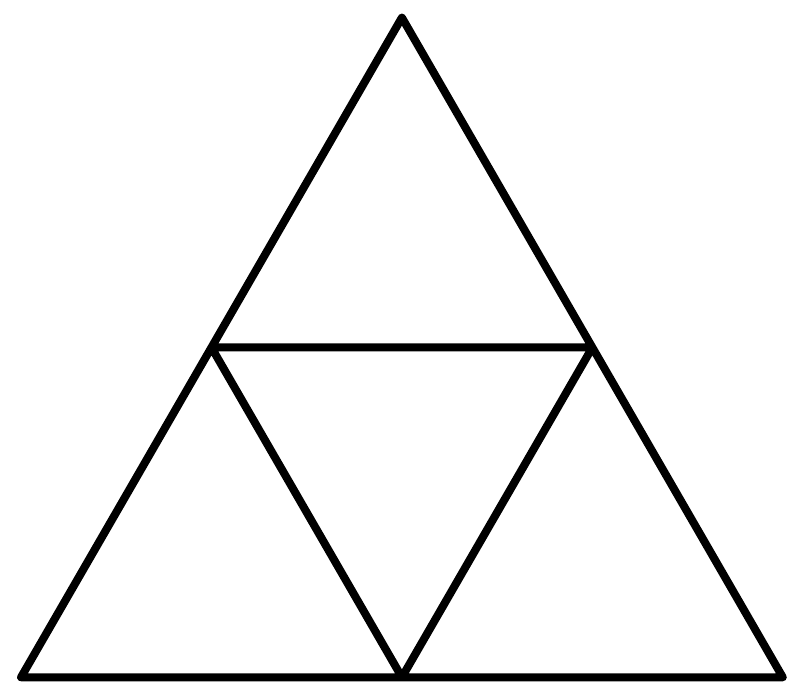}
\caption{A triangular foldable figure.}\label{f1000}
\end{center}
\end{figure}

These observations suggest now the validity of the following more general situation.
\begin{proposition}
Suppose that $C$ is any foldable figure which is a fundamental domain for both a translation group $\mathcal{T}$ and the affine Weyl group $\widetilde{\mathcal{W}}$ for $C$, and which contains 0 in its interior. If the intersection group $\mathcal{J}$ of $\mathcal{T}$ and $\widetilde{\mathcal{W}}$ is big enough, in the sense described above, so that for an expansive matrix $A\in M_n (\mathbb{R})$ the dilation group $\mathcal{D}$ for $A$ and the intersection group $\mathcal{J}$ is an abstract dilation--translation pair, then there exist sets $W$ which are simultaneously dilation-translation and dilation--reflection wavelet sets.
\end{proposition}
\begin{proof}
Let $B$ denote the unit ball of $\R^n$ and let $W$ be the set whose existence is guaranteed by Theorem \ref{t1} (Theorem 1 of  \cite{DLS1}) and which is both dilation--congruent to $F := A(B) \setminus B$  and $\mathcal{J}$-congruent to $C$. (Here $B$ is the unit ball of $\mathbb{R}^n$.) Then $W$ is automatically congruent to $C$ by the larger groups, namely the affine Weyl group $\widetilde{\mathcal{W}}$ and the translation group $\mathcal{T}$. Now, dilation--congruence to $F$ and translation-congruence to $C$ makes $W$ are dilation--translation wavelet set, whereas dilation--congruence to $F$ and $\widetilde{\mathcal{W}}$--congruence to $C$ make it a dilation--reflection wavelet set.
\par
The key is that $W$ needs to tile $\mathbb{R}^n$ under the full translation group $\mathcal{T}$, and also the full affine Weyl group $\widetilde{\mathcal{W}}$. Since this is true for $C$, the congruences of $C$ to $W$ under both groups guarantees that $W$ tiles $\mathbb{R}^n$, too.
\end{proof}
\vskip 5pt
The two wavelets presented in Examples \ref{ex1} and \ref{example2} are very special, for the reasons given above. In addition, the affine Weyl group associated with the foldable figure $C = [-\pi,\pi]\times [-\pi,\pi]$ is also rather specific in the sense that it is reducible. It consists of the bifold product of the Weyl group associated with the interval $[-\pi,\pi]$. Theorem \ref{100} guarantees the existence of dilation-reflection wavelet sets in general, even when the Weyl group is irreducible.  So it may be interesting to construct concrete examples of such sets. One example of a foldable figure whose affine Weyl group is irreducible is given in Figure \ref{f1000}.  We therefore pose the following problem.
\vskip 5pt\noindent
\textsc{PROBLEM 1:} Given the foldable figure $C$ depicted in Figure \ref{f1000}, construct concrete examples of dilation--reflection wavelet sets $W$ for $C$. In particular, are there any such examples which are bounded and bounded away
from 0 (as are the sets in Examples 7.7 and 7.8). In principle, if one follows the proof of Theorem 2.6 in [DLS1] using that constructive proof as an algorithm for constructing a wavelet set, then dilation-reflection sets for the Weyl group associated with Figure \ref{f1000} can easily be constructed.  However, with that method, the sets constructed are not bounded subsets of the plane, and are also not bounded away from 0. In addition, they are difficult to work with.
\par
Another interesting question connecting with the above arguments is the following.
\vskip 5pt\noindent
\textsc{PROBLEM 2:} Let $C$ be any foldable figure in $\mathbb{R}^n$ containing 0 in its interior and let $\widetilde{\mathcal{W}} = \mathcal{W}\ltimes \Gamma$ be the associated affine Weyl group. Suppose that $A$ is any expansive matrix in $M_n (\mathbb{R})$, and $\mathcal{T}$ a translation group on $\mathbb{R}^n$ obtained by translating with respect to a basis $\mathcal{B} := \{b_1, b_2,...b_n\}$ for $\mathbb{R}^n$, i.e., a vector $x$ in $\mathbb{R}^n$ is mapped to $x - (k_1b_1,  k_2b_2, ..., k_nb_n)$ for all $n$-tuples $(k_1, ..., k_n)\in\mathbb{Z}^n$. Give necessary and sufficient conditions for the existence a set $W$ which is simultaneously
\begin{enumerate}
\item $\widetilde{\mathcal{W}}$--congruent to $C$;
\item $D_A $--congruent to the set $F := A(B) \setminus B$  (as in the proof of Theorem \ref{100});
\item $\mathcal{T}$--congruent to the set $[0, b_1)\times [0, b_2)\times \ldots\times[0, b_n)$.
(Note that this last set is the simplest fundamental domain for $\mathcal{T}$).
\end{enumerate}
Any $W$ satisfying (1), (2), and (3) would be both a dilation--translation wavelet set for $(\mathcal{D},\mathcal{T})$ and a dilation--reflection wavelet set for $(\widetilde{\mathcal{W}}, A)$.  Conversely, any set which is both a dilation-translation wavelet set for $(\mathcal{D},\mathcal{T})$ and a dilation--reflection wavelet set for $(\widetilde{\mathcal{W}}, A)$ must satisfy (1), (2), and (3). In particular, does there exist such a $W$ for an irreducible Weyl group, such as the group corresponding to the foldable figure in Figure \ref{f1000}?  We think that the answer is probably no.  But in the topic of wavelet sets there are often surprises, so we would not be very surprised if the answer was yes.
\section*{Acknowledgments}\noindent
The authors thank the anonymous referee for their detailed comments and recommendations improving the understanding of the paper and some interesting suggestions for further work.  In particular, the referee observed that because of the semi-direct product structure of the affine Weyl group and the existence of three-way tiling sets, there may be connections between the dilation-reflection wavelet theory developed in this paper and the composite dilations developed by Krishtal, Robinson, Weiss and Wilson in \cite{KRWW}.

\bibliographystyle{amsalpha}

\end{document}